\documentclass{statsoc}
\usepackage{graphicx,subfigure}
\usepackage{amssymb}
\usepackage{color}

\usepackage{geometry}
\usepackage{amsmath}
\usepackage{amsbsy}
\usepackage{enumerate}
\usepackage{bm}

\newtheorem{assumption}{Assumption}
\newtheorem{example}{Example}

\newtheorem{corollary}{Corollary}

\newtheorem{proposition}{Proposition}

\newtheorem{defi}{Definition}
\newtheorem{thm}{Theorem}
\newtheorem{lemma}{Lemma}

\def\ind{\begin{picture}(9,8)
         \put(0,0){\line(1,0){9}}
         \put(3,0){\line(0,1){8}}
         \put(6,0){\line(0,1){8}}
         \end{picture}
        }

\title[Principal stratification estimation and surrogate endpoint evaluation]{Principal causal effect identification and surrogate endpoint evaluation by multiple trials}
\author{Zhichao Jiang}
\address{Peking University, Beijing, People's Republic of China}
\author{Peng Ding}
\address{Harvard University,
         Cambridge, Massachusetts,
         USA.}
\author[Jiang {\it et al.}]{Zhi Geng}
\address{Peking University, Beijing, People's Republic of China}

\begin{document} 

\footnotetext{{\\ \it Address for correspondence}: Peng Ding, Department of Statistics, Harvard University, One Oxford Street, Cambridge 02138 Massachusetts, U.S.A.\\
E-mail: pengdingpku@gmail.com}
         
\begin{abstract}
Principal stratification is a causal framework to analyze randomized experiments with a post-treatment variable  between the treatment and endpoint variables. Because the principal strata defined by the potential outcomes of the post-treatment variable are not observable, we generally cannot identify the causal effects within principal strata. Motivated by a real data set of phase III adjuvant colon clinical trials, we propose approaches to identifying and estimating the principal causal effects via multiple trials. For the identifiability, we remove the commonly-used exclusion restriction assumption by stipulating that the principal causal effects are homogeneous across these trials. To remove another commonly-used monotonicity assumption, we give a necessary condition for the local identifiability, which requires at least three trials. Applying our approaches to the data from adjuvant colon clinical trials, we find that the commonly-used monotonicity assumption is untenable, and disease-free survival with three-year follow-up is a valid surrogate endpoint for overall survival with five-year follow-up, which satisfies both the causal necessity and the causal sufficiency. We also propose a sensitivity analysis approach based on Bayesian hierarchical models to investigate the impact of the deviation from the homogeneity assumption. 
\end{abstract}


\section{Introduction}

\subsection{Principal stratification and surrogate endpoints}

Causal effects are defined as comparisons between potential outcomes 
of the endpoint variable under treatment and control for the same group of individuals.
To evaluate the causal effects on the endpoint within subpopulations
stratified by a post-treatment variable, \citet{frangakis2002principal} proposed the principal stratification framework
in which the subpopulations are defined by the joint values of the potential outcomes 
instead of the observed value of the post-treatment variable.
The joint potential outcomes of the post-treatment variable 
can be viewed as a pretreatment covariate vector unaffected by the treatment.
Principal stratification has several applications in the current literature.
The compliance behavior defined by the potential treatment
acceptances is used to address the noncompliance problem \citep{angrist1996identification}, 
the potential survival status is used to evaluate the effect on the quality of life with truncation by death (\citealp{rubin2006causal}; \citealp{ding2011identifiability}; \citealp{yang2015using}), and
the potential employment status is used to deal with the truncation of wages due to unemployment in  the evaluation of job-training programs \citep{zhang2003estimation, zhang2009likelihood, frumento2012evaluating}.  The potential response indicators are used to address non-ignorable nonresponse problems \citep{frangakis1999addressing,  mealli2008comparing, mattei2014identification}, and the potential intermediate variables are used to define direct and indirect effects \citep{rubin2004direct,  gallop2009mediation, mattei2011augmented}.

Another important application of principal stratification is to evaluate surrogate endpoints in clinical trials \citep{frangakis2002principal}.
When direct measurement of an endpoint of interest is too time-consuming or costly, we try to measure a surrogate for the endpoint. There have been some criteria for judging surrogates from different points of view. \citet{prentice1989surrogate} first proposed the statistical surrogate criterion, which requires conditional independence of the treatment and the observed endpoint given the observed surrogate. From a causal perspective,
\citet{frangakis2002principal} proposed the principal surrogate criterion, and 
pointed out that the statistical surrogate
does not satisfy the causal necessity, i.e., no difference between the potential outcomes of the surrogate implies no difference between the potential outcomes of the endpoint. \citet{lauritzen2004discussion} proposed the strong surrogate criterion depicted
by a causal diagram requiring the surrogate break the causal path
from the treatment to the endpoint, which is stronger than the principal surrogate. \citet{gilbert2008evaluating} argued that a principal surrogate should satisfy not only causal necessity but also causal sufficiency. The causal sufficiency requires that if the treatment effect on the surrogate is non-zero, then the treatment effect on the endpoint is also non-zero. \citet{joffe2009related} considered four different approaches for evaluating surrogate endpoints. However, all of these approaches may suffer from the surrogate paradox as pointed out by \citet{chen2007criteria}, \citet{ju2010criteria} and \citet{vanderweele2013surrogate}.
That is, for a surrogate satisfying any of these criteria, it is possible that the treatment has a negative average causal effect on the endpoint
even if the treatment has a positive average causal effect on the surrogate and the surrogate has a positive average causal effect on the endpoint.

For evaluating surrogates, \citet{mealli2012refreshing} suggested conducting a principal stratification analysis investigating the effect of the treatment within all principal strata.
However, because we cannot simultaneously observe two potential outcomes
of the post-treatment variable,
we do not know the principal stratum of an individual, and generally we cannot identify causal effects within principal strata.
\citet{zhang2003estimation} and \citet{cheng2006bounds} proposed large sample bounds of causal effects within some principal strata, which, however, may be too wide to be informative. \citet{angrist1996identification} discussed the identifiability of the complier average causal effect under the monotonicity and exclusion restriction assumptions. Without the exclusion restriction assumption, \citet{zhang2009likelihood} used Gaussian mixture models to identify causal effects within principal strata.
\citet{gilbert2008evaluating} and \citet{huang2011comparing}  proposed approaches to evaluating surrogates based on principal stratification in a single trial, but they assumed constant potential outcomes of the surrogate under control. \citet{zigler2012bayesian} proposed a Bayesian approach to estimating the causal effects on the endpoint within principal strata without the monotonicity assumption, but their approach relied on prior distributions on some parameters that are not identifiable.

\subsection{Motivation}

Our study is motivated by the data from phase III adjuvant colon clinical trials (ACCTs).
The goal of the ACCTs is to test whether disease-free survival (DFS) with three-year follow-up can be used
as a surrogate for the overall survival (OS) with five-year follow-up \citep{sargent2005disease}.
The data are collected from 10 clinical trials. \citet{sargent2005disease} found a strong correlation between the hazard ratio of the treatment on DFS and the hazard ratio of the treatment on OS.
Applying a meta-analysis, \citet{baker2012predicting} used DFS as a principal surrogate to predict the treatment effect on the endpoint using the treatment effect on the surrogate. They assumed that the true endpoint does not depend on the treatment among subjects whose surrogate is the same under treatment and control, i.e., the causal necessity holds.
However, the causal necessity used by \citet{baker2012predicting} has not been verified by either previous studies or observed data.

For the ACCTs data,
the monotonicity assumption requiring nonnegative individual causal effect on the surrogate may not hold, because the treatment may have negative side-effects on the surrogate for some patients.
Second, the exclusion restriction assumption implies the causal necessity, which is the scientific question of interest in evaluating surrogate endpoints.
Third, Gaussian mixture models are not applicable to the binary endpoints in the ACCTs data. The meta-analysis approaches \citep{daniels1997meta,li2011causal,baker2012predicting} may also suffer from the surrogate paradox \citep{vanderweele2013surrogate}. 
When the monotonicity or exclusion restriction assumption does not hold, \citet{ding2011identifiability},
\citet{mealli2013using}, and \citet{yang2015using} achieved partial and point identification of principal causal effects by exploiting covariates or secondary outcomes, and \citet{mattei2012exploiting} improved Bayesian inference for principal causal effects using multiple outcomes. However, when additional outcomes or covariates are not available, we cannot identify causal effects using previous approaches. Without identifiability, Bayesian inference is sensitive to  prior distributions \citep{gustafson2009limits}.

In this paper, we propose approaches to identifying the principal stratification causal effects without the exclusion restriction assumption
and further without the monotonicity assumption.
To remove the exclusion restriction assumption,
we need at least two trials.
Furthermore, to remove the monotonicity assumption,
we give a sufficient condition for the local identifiability 
of the causal effects, which requires at least three trials.
We first assume that these trials are homogeneous.
Using the identified principal stratification causal effects,
we can test the causal necessity and the causal sufficiency
of a surrogate.
We evaluate the treatment effect
on the unobserved endpoint in a new trial where the distribution of principal strata
and/or the treatment may be different
from those in the validation trials. We then propose criteria for surrogates that avoid the surrogate paradox.
Allowing for possible deviations from the homogeneity assumption, we then conduct a sensitivity analysis based on a class of Bayesian hierarchical models.
We apply the proposed approaches to the ACCTs data
and evaluate the surrogacy of DFS with three-year follow-up
for OS with five-year follow-up.

The paper proceeds as follows.
We introduce the notation and assumptions in Section \ref{sec::notation-assumptions-bounds}. We present the identification conditions with and without the monotonicity assumption in Section \ref{sec::identification-multiple-trials}.
Section \ref{sec::evalution-principal-surrogacy} discusses the evaluation of surrogate based on principal stratification.
We evaluate the performance of
our approaches with finite sample sizes in Section \ref{sec::simulation-studies} via simulation studies.
We apply our approaches to the ACCTs data in Section \ref{sec::application-ACCT}, and perform a sensitivity analysis
using Bayesian hierarchical models in Section \ref{sec::sensitivity}.
We conclude with a discussion in Section \ref{sec::discussion}. We present the
details of the proofs and computations in the online supporting materials, and provide the data and compute code online.

\section{Notation and assumptions}
\label{sec::notation-assumptions-bounds}

\subsection{Potential outcomes and principal stratification}

Suppose that there are $N_R$ independent trials
which have the same treatment and the same control.
Let $R$ denote the trial number taking values from $1$ to $N_R$.
Let $Z$ be a binary treatment assignment with
$Z=1$ for treatment and 0 for control.
Let $S$ denote a binary surrogate
and $Y$ a binary endpoint of interest.
In the ACCTs data,
$S_i = 0$ if the cancer of patient $i$ reoccurred before $3$ years, and $S_i=1$ otherwise.
The endpoint $Y_i$  denotes the survival status within $5$ years, with $Y_i = 1$ for ``survival'' and $Y_i = 0$ for ``dead.''
We use the potential outcomes framework to define causal effects and make
the Stable Unit Treatment Value Assumption (SUTVA) throughout the paper, i.e., there is only one version of the potential outcomes and there is no interference between units \citep{rubin1980randomization}. The SUTVA allows us to uniquely define the potential outcomes $S_i(z)$ and $Y_i(z)$ for the surrogate and the endpoint variables if patient $i$ were to receive treatment $z$ for $z=0$ and $1$.
The observed values of $S$ and $Y$ are deterministic functions of both the treatment assignment and the potential outcomes, i.e., $S_i = S_i(Z_i) = Z_iS_i(1) + (1-Z_i)S_i(0)$ and $Y_i = Y_i(Z_i) = Z_iY_i(1) + (1 - Z_i)Y_i(0)$.
Throughout our paper, we assume that $\{ (R_i, Z_i, S_i(1), S_i(0), Y_i(1), Y_i(0))  : i=1, \ldots, N \}$ are independently and identically distributed (iid) samples drawn from an infinite super-population, and thus the observed $\{ (R_i ,Z_i, S_i, Y_i) : i=1,\ldots, N  \}$ are also iid.
Although we will discuss only the binary endpoints in depth, our framework and approaches are also applicable to general continuous endpoints by dichotomizing the endpoints to identify the distributional causal effects.

\citet{frangakis2002principal} defined principal stratification using the joint potential outcomes of the surrogate under both treatment and control, i.e., $U_i = (S_i(1), S_i(0))$. For simplicity, we relabel the possible values of $U$,
 $(1,1), (1,0), (0,1)$ and $(0,0)$,  as $ss, s\bar{s}, \bar{s}s$ and $\bar{s}\bar{s}$, respectively.
In our motivation example, we use ``$s$'' for ``disease-free survival in three-year follow-up,'' and ``$\bar{s}$'' otherwise.
We are interested in the principal stratification average causal effects (PSACEs) for each trial, defined as
$$
ACE_{ur} = E\{ Y(1)-Y(0)\mid  U=u, R=r\}
$$
for $u=ss, s\bar{s}, \bar{s}s, \bar{s}\bar{s}$, and $r=1, \ldots, N_R$.
Based on the PSACEs, we can formally define the causal necessity and causal sufficiency (Gilbert and Hudgens 2008).
\begin{defi}
Causal necessity requires $ACE_{ur}=0$ for $u=ss,\bar{s}\bar{s}$, and
causal sufficiency requires $ACE_{ur}\neq 0$ for $u=s\bar{s},\bar{s}s$, where $r=1, \ldots, N_R$.
\end{defi}

The above definition of causal necessity is weaker than the usual exclusion restriction assumption \citep{angrist1996identification} that requires zero individual causal effect on the outcome for principal strata $u=ss$ and $\bar{s}\bar{s}.$

\subsection{Assumptions}
\label{sec::assumption}

In this subsection, we introduce the basic assumptions, which are commonly used in causal inference.
By identification of parameters, we mean that they can be expressed as functions of the distributions of observed variables. We let $A\ind B \mid C$ denote the conditional independence of $A$ and $B$ given $C$.

\begin{assumption}
\label{assume:2}
(Randomization). $Z \ind \{S(1),S(0),Y(1),Y(0) \}  \mid R$.
\end{assumption}

Assumption \ref{assume:2} means that treatment assignment $Z$ is randomized
within each trial, but the assignment probabilities may be different.
Because we have independent randomized clinical trials in the ACCTs data, Assumption \ref{assume:2} holds by the designs of experiments.

The following monotonicity assumption is widely used for
principal stratification analysis \citep{angrist1996identification} and principal surrogate evaluation \citep{gilbert2008evaluating}.

\begin{assumption}
\label{assume:3}
(Monotonicity). $S_i(1)\geq S_i(0)$ for each individual $i$.
\end{assumption}

Assumption \ref{assume:3} means that the treated surrogate endpoint $S_i(1)$ is always better than or equal to
the controlled surrogate endpoint $S_i(0)$ for every patient $i$.
The monotonicity assumption may be too restrictive for some real applications where treatment may have negative side-effects on some patients.
We shall discuss the identification of PSACEs with and without the monotonicity assumption separately.

Let $O(z,s)$ denote the set of principal strata that are
compatible with the observed values $Z_i=z$ and $S_i=s$.
Then we have $O(1,1) = \{ss, s\bar{s}\}$, $O(1,0) = \{ \bar{s}\bar{s}, \bar{s}s\}$, $O(0,1) = \{ ss, \bar{s}s\}$, and $O(0,0) = \{ s\bar{s}, \bar{s}\bar{s}\}$ without the monotonicity assumption.
The monotonicity assumption
eliminates the stratum $\bar{s}s$, and therefore $O(1,1) = \{ss, s\bar{s}\}$, $O(1,0) = \{ \bar{s}\bar{s}\}$, $O(0,1) = \{ ss\}$, and $O(0,0) = \{ s\bar{s}, \bar{s}\bar{s}\}$.
Although the principal stratification variable $U$ is not observable for all units, monotonicity allows us to
identify the proportions of principal strata using the distribution of the observed data. Define 
$$
\pi_{ur}=P(U=u\mid R=r),
$$
which can be identified by
$\pi_{\bar{s}\bar{s},r} = P(S=0\mid Z=1,R=r)$,
$\pi_{ss,r} = P(S=1\mid Z=0,R=r)$
and $\pi_{s\bar{s},r}  = 1-\pi_{ss,r}-\pi_{\bar{s}\bar{s},r}$.
With monotonicity, we can also identify the expectations
$E\{Y(0)\mid U=ss,R=r\} = E(Y \mid Z=0,S=1,R=r)$ and $E\{Y(1)\mid U=\bar{s}\bar{s},R=r\} = E(Y\mid Z=1,S=0,R=r)$.
Since patients within the observed group $(Z=1, S=1)$ or $(Z=0, S=0)$ are both mixtures of two latent principal strata, we cannot identify the PSACEs without further assumptions beyond Assumptions \ref{assume:2} and \ref{assume:3}. We can obtain large sample bounds  
for the PSACEs from the observed data, but these bounds are barely informative as shown in the online Appendix.

\section{Identification of the PSACEs from multiple trials}
\label{sec::identification-multiple-trials}

Because the principal stratum is unobservable,
causal effects within principal strata
are not identifiable in general.
In this section, we shall propose approaches
to improving the identifiability of the PSACEs
in terms of multiple trials.
To combine information from multiple trials,
we  make the following homogeneity assumption.

\begin{assumption}
\label{assume:4}
(Homogeneity). $R \ind Y(z)  \mid U$ for $z=0,1$.
\end{assumption}

The homogeneity assumption means
that the potential outcome $Y(z)$ may depend on
the principal stratum $U$,
but is independent of the trial number $R$
conditional on the principal stratum $U$. Thus, there is no ``direct effect'' of $R$ on $Y(z)$ within the principal stratum $U$.  We can include pretreatment covariates $X$  to make this assumption more plausible, i.e., $R \ind Y(z)  \mid (U,X)$, and we omit it for simplicity. As an example,  Assumption \ref{assume:4} for the ACCTs data  means that the survival of a patient under treatment or control does not depend on which trial he/she is in once we know  his/her principal stratum $U$, defined by the reoccurrence of cancer under both treatment and control.  In particular, for the patients in stratum $U = ss$, the cancer will not reoccur within 3 years
no matter whether they receive the treatment or not.  Thus the principal stratum can be interpreted as a measure of the physical status of a patient. It may be plausible that  the survival will no longer depend on the trial number  after we know a patient's physical status and the cancer status. Recognizing that Assumption \ref{assume:4} may not be directly testable and may be violated in practice, we propose an approach for conducting sensitivity analysis  
about Assumption \ref{assume:4} in Section \ref{sec::sensitivity}.

The following result further explains the homogeneity of causal effects
across trials.

\begin{proposition}
\label{prop:2}
Under Assumption \ref{assume:2}, Assumption \ref{assume:4} is equivalent to $R \ind Y  \mid (U, Z)$, which implies $ACE_{ur} = ACE_{ur'}$ for all $u=ss$,
$s\bar{s}$, $\bar{s}s$, $\bar{s}\bar{s}$ and $r \neq r'$.
\end{proposition}

 By Proposition \ref{prop:2}, we can write $ACE_{u} = ACE_{ur}$,
and thus we can discuss the identifiability of the common PSACEs
using the data from multiple trials together.

To simplify the notation,
we define $p_r = P(R=r)$ as the proportion of trial $r$, $\alpha_r = P(Z=1\mid R=r)$ as the proportion of patients receiving treatment in trial $r$,  $\delta_{zur} = P(Y=1\mid Z=z, U=u, R=r)$ as the survival proportion conditional on the treatment $z$, principal stratum $u$ and trial number $r$. Under the homogeneity assumption, we have that $\delta_{zu}=\delta_{zur}$. 
We further define $\bm{p} = \{p_r: r=1,\ldots,N_R\}$, $\bm{\alpha}=\{\alpha_r: r=1,\ldots,N_R\}$,
$\bm{\pi}_r=\{\pi_{ss,r}, \pi_{s\bar{s},r}, \pi_{\bar{s}\bar{s},r}, \pi_{\bar{s}s,r} \}$,
$\bm{\pi}=\{  \bm{\pi}_{r}: r=1,\ldots,N_R\}$,
and $ \bm{\delta}=\{\delta_{zu }: z=0,1; u=ss,s\bar{s},\bar{s}\bar{s}, \bar{s}s  \}$.
Under monotonicity, we do not have the parameters $\delta_{z,\bar{s}s}$, and we have $\pi_{\bar{s}s,r}=0$.

\subsection{\bf Identification with monotonicity}

In addition to the homogeneity assumption,
we need the following assumption.

\begin{assumption}
\label{assume:5}
(a) There exist at least two trials $r_1$ and $r_2$, such that $\pi_{ss,r_1} / \pi_{s\bar{s},r_1} \neq \pi_{ss,r_2} / \pi_{s\bar{s},r_2}$.

(b) There exist at least two trials $r_3$ and $r_4$, such that $\pi_{s\bar{s},r_3} / \pi_{\bar{s}\bar{s},r_3} \neq \pi_{s\bar{s},r_4} / \pi_{\bar{s}\bar{s},r_4}$.
\end{assumption}

Again, treating the principal stratum as a measure of patients' physical status, 
Assumption \ref{assume:5} means that patients in different trials have different distributions of the patients' physical status. Under monotonicity, Assumption \ref{assume:5} is testable by data because we can identify the proportions of the principal strata under Assumption \ref{assume:3} as shown in Section 2.2.

The trial number $R$ acts like an instrumental variable associated with $U$ in the sense that Assumption 3 is similar to the exclusion restriction assumption, and Assumption 4 guarantees the association between $R$ and $U$. 
Under Assumptions 3 and 4, there is no ``direct effect'' of $R$ on $Y$, which is similar to the instrumental variable case where the instrumental variable does not ``directly'' affect the outcome.

\begin{thm}
\label{thm:1}
Under Assumptions \ref{assume:2} to \ref{assume:4}, we have that, for $u=ss$, $s\bar{s}$ and $\bar{s}\bar{s}$,
\begin{enumerate}
[(a)]
\item
$P(Y=1\mid Z=1,U=u)$ for all principal strata $u$ are identifiable if Assumption \ref{assume:5}(a) holds;

\item
$P(Y=1\mid Z=0,U=u)$ for all principal strata $u$ are identifiable if Assumption \ref{assume:5}(b) holds; 

\item
$ACE_{u}$ for all principal strata $u$ are identifiable if both Assumptions \ref{assume:5}(a) and \ref{assume:5}(b) hold.
\end{enumerate}
\end{thm}

Assumption \ref{assume:5} requires at least two independent trials
for the identifiability. We can use any two trials satisfying Assumption \ref{assume:5} to obtain moment estimators, and we can also use all trials to obtain more efficient maximum likelihood estimates (MLEs).
The following example illustrates
the identifiability for the case with two independent trials $(N_R = 2)$ under Assumptions \ref{assume:2} to \ref{assume:5}.

\begin{example}
First, the proportions of the principal strata $\pi_{ur}$
are identifiable under monotonicity, as discussed in Section 2.2.

 Second, the probability $\omega_{ys \mid zr} = P(Y=y, S=s\mid Z=z, R=r)$ can be identified from the observed data. We can directly identify two outcome distributions
$\delta_{1,\bar{s}\bar{s}}=  \omega_{10 \mid 11} / \pi_{\bar{s}\bar{s},1}$
and
$\delta_{0,ss}= \omega_{11 \mid 01} / \pi_{ss,1}$.

Third, if $\pi_{ss,1}/\pi_{s\bar{s}, 1} \neq \pi_{ss, 2}/\pi_{s\bar{s}, 2}$, we have from the proof of Theorem \ref{thm:1} that
\begin{eqnarray*}
\delta_{1,ss} = \frac{\omega_{11\mid 11}\cdot \pi_{s\bar{s},2}-\omega_{11\mid 12}\cdot \pi_{s\bar{s},1}}{\pi_{ss,1}\cdot \pi_{s\bar{s},2}-\pi_{ss,2}\cdot \pi_{s\bar{s},1}}, \quad
\delta_{1,s\bar{s}} = \frac{\omega_{11\mid 12}\cdot \pi_{ss,1}-\omega_{11\mid 11}\cdot \pi_{ss,2}}{\pi_{ss,1}\cdot \pi_{s\bar{s},2}-\pi_{ss,2}\cdot \pi_{s\bar{s},1}}.
\end{eqnarray*}
If $\pi_{ss,1}/\pi_{s\bar{s}, 1} \neq \pi_{ss, 2}/\pi_{s\bar{s}, 2}$, we have from the proof of Theorem \ref{thm:1} that
\begin{eqnarray*}
\delta_{0,s\bar{s}} = \frac{\omega_{10\mid 01}\cdot \pi_{\bar{s}\bar{s},2}-\omega_{10\mid 02}\cdot \pi_{\bar{s}\bar{s},1}}{\pi_{s\bar{s},1}\cdot \pi_{\bar{s}\bar{s},2}-\pi_{s\bar{s},2}\cdot \pi_{\bar{s}\bar{s},1}}, \quad
\delta_{0,\bar{s}\bar{s}} = \frac{\omega_{10\mid 02}\cdot \pi_{s\bar{s},1}-\omega_{10\mid 01}\cdot \pi_{s\bar{s},2}}{\pi_{s\bar{s},1}\cdot \pi_{\bar{s}\bar{s},2}-\pi_{s\bar{s},2}\cdot \pi_{\bar{s}\bar{s},1}}.
\end{eqnarray*}
\end{example}

\subsection{\bf Identification without monotonicity}

Monotonicity for every individual may be too restrictive in practice.
In this subsection, with the help of more than two independent trials, we can remove monotonicity.

By the homogeneity assumption, the probabilities of the observed data can be decomposed as
\begin{eqnarray}\label{eq::obs}
&&P(Z = z, S=s, Y=y\mid R=r) \nonumber \\ 
&=& \sum_{u\in O(z,s)}  P(Z=z \mid R=r) \cdot   \pi_{ur} \cdot  P(Y=y\mid Z=z, U=u),
\end{eqnarray}
for $z, s, y = 0,1$.
Let $\boldsymbol{\xi}=( \bm{p}, \bm{\alpha},  \bm{\pi}, \bm{\delta})$ denote the vector of the parameters.
Let $\bm{f}$ be the vector of probabilities on the left-hand side of equation (\ref{eq::obs}).
Applying a Taylor expansion, we can approximate $\bm{f} (\boldsymbol{\xi}) $ by linear equations of $\boldsymbol{\xi}$ around the true parameters $\boldsymbol{\xi_0}$, i.e.,
$
 \boldsymbol{f} (\boldsymbol{\xi})  \approx \boldsymbol{f} ( \boldsymbol{\xi_0} ) +\nabla \boldsymbol{f} \mid _{\boldsymbol{\xi_0}} (\boldsymbol{\xi}-\boldsymbol{\xi_0}),
$
where $\nabla \boldsymbol{f} \mid _{\boldsymbol{\xi_0}}$ is the Jacobian matrix with the $(i,j)$-th element $\partial f_i/ \partial \xi_j \mid _{\boldsymbol{\xi_0}}$.
According to ~\citet{bandeen1997latent},
a distribution $F_W(w;\boldsymbol{\xi})$ of a random variable $W$ is locally identifiable at $\boldsymbol{\xi_0}$ if there exists some neighborhood $N(\boldsymbol{\xi_0})$ of $\boldsymbol{\xi_0}$ such that for all $ \boldsymbol{\xi} \in N(\boldsymbol{\xi_0})$, $F_W(w;\boldsymbol{\xi_0})=F_W(w;\boldsymbol{\xi})$ for all $w$ if and only if $\bm{\xi}=\bm{\xi}_0$. Therefore,
under the randomization and homogeneity assumptions,
the parameter vector $\boldsymbol{\xi}$ can be represented by
a function of the distributions of observed variables
and thus it is locally identifiable if the Jacobian matrix $\nabla \boldsymbol{f} \mid _{\boldsymbol{\xi_0}}$ is of full column rank. In fact, Assumption \ref{assume:5} is a necessary condition for $\nabla \boldsymbol{f} \mid _{\boldsymbol{\xi_0}}$ to be full column rank. In practice,
the full rank condition of $\nabla \boldsymbol{f} \mid _{\boldsymbol{\xi_0}}$
can be tested empirically by the rank of $\nabla \boldsymbol{f} \mid _{\boldsymbol{\hat{\xi}}}$, where $\boldsymbol{\hat{\xi}}$ is the MLE of $\boldsymbol{\xi_0}$ \citep{goodman1974exploratory, skrondal2004generalized}.

A necessary condition for local identifiability of
unknown parameters is that the number of observed frequencies is larger than the number of unknown parameters. For this case, we have the following result.

\begin{proposition}
\label{prop:nec}
Under Assumptions \ref{assume:2} and \ref{assume:4},
a necessary condition for the local identifiability of the joint distribution of  $(Y,Z,S,U,R)$ is $N_R\geq 3$.
\end{proposition}

Intuitively, when $N_R\leq  2$, there will be less equations than unknown parameters in the equations in (\ref{eq::obs}), and thus we cannot obtain unique solution of the unknown parameters. Therefore, there must be at least three independent trials to identify the PSACEs without monotonicity.
In the case without monotonicity, we can not obtain closed forms of PSACEs  like Example 1, because we are not able to obtain closed forms of $\pi_{ur}$'s. We can use a numerical approach to calculate the PSACEs, and we give an example in the online Appendix to illustrate the identifiability for the case without monotonicity.
 
Without monotonicity, we can prove only local identification. If a parameter is locally identifiable but not globally identifiable, its posterior distribution must be multimodal and has the same value at multiple modes. Thus, we can check its posterior distribution to verify identifiability. In our simulation studies and application, posterior distributions of all the parameters are  unimodal, which means that the parameters are indeed globally identifiable. In general, we can test whether the posterior distributions of parameters are unimodal using existing methods.

\subsection{\bf Computation, model checking, and goodness-of-fit tests}
\label{sec::computation}

In this subsection,
we shall discuss the estimation of the parameters.
We can use the expectation-maximization (EM) algorithm and the Gibbs Sampler to compute the MLEs and simulate the posterior distributions of PSACEs.
Suppose that $R$ follows a categorical distribution with parameter $\bm{p}$, $Z$ follows Bernoulli distributions conditional on $R$ with parameter $\bm{\alpha}$, $U$ follows categorical distributions conditional on $R$ with parameter $\bm{\pi}$, and $Y$ follows Bernoulli distributions conditional on $Z$ and $U$ with parameter $\bm{\delta}$.
The unobserved variable $U$ is treated as a latent variable.
The complete data can be represented by a contingency table classified by
$(Z, U, Y, R)$, and the observed data can be represented by a contingency table classified by $(Z, S, Y, R)$ with cell counts $N_{zsyr}=\#\{i: Z_i=z,S_i=s,Y_i=y,R_i=r \}.$ 
For Bayesian inference,
we use Dirichlet (or Beta) distributions
with parameters $(1, \ldots, 1)$
as the non-informative prior distributions
of $\bm{p}$, $\bm{\alpha}$, $\bm{\pi}$ and $\bm{\delta}$.
Computational details are given in the online Appendix.

More interestingly, when $N_R\geq 2$ with monotonicity or $N_R\geq 3$ for without monotonicity, the extra degrees of freedom allow us to perform goodness-of-fit tests based on the asymptotic distributions of the likelihood ratio statistics:
\begin{eqnarray*}
&& T_m(N_{zsyr})  = 2\log   \frac{L(\boldsymbol{\hat{\zeta}}_s \mid N_{zsyr})}{L(\boldsymbol{\hat{\xi}}_{m}\mid N_{zsyr})}
\stackrel{a}{\sim} \chi^2_{  4N_R - 6 } , \quad \\
&& T_{nm}(N_{zsyr})  =  2\log   \frac{L(\boldsymbol{\hat{\zeta}}_s \mid N_{zsyr})}{L(\boldsymbol{\hat{\xi}}_{nm}\mid N_{zsyr})}
\stackrel{a}{\sim} \chi^2_{  3N_R - 8  },
\end{eqnarray*}
where $\boldsymbol{\hat{\zeta} } _s$ is the MLE under the saturated model, i.e., $(Z,S,R,Y)$ follows the multinomial distribution without constraints on parameters, and $\boldsymbol{\hat{\xi}}_m$ and $\boldsymbol{\hat{\xi}}_{nm} $ are the MLEs with and without monotonicity.
The degrees of freedom with and without monotonicity are $ (8N_R - 1) - (4N_R+5) = (4N_R-6)$ and $  (8N_R - 1) -  (5N_R + 7)=  (3N_R - 8)$, respectively.
On the other hand, we can use Bayesian posterior predictive $p$-values (ppp)
for model checking \citep{rubin1984bayesianly, meng1994posterior}.
Let $N_{zsyr}^{\text{rep}}$ denote the replicate data generated from their posterior predictive distributions.
The ppp's for the models with and without monotonicity are
\begin{eqnarray*}
&& ppp_m= P\{T_m(N_{zsyr})>T_m(N_{zsyr}^{\text{rep}})\mid  N_{zsyr}\},\\
&& ppp_{nm} = P\{T_{nm}(N_{zsyr})>T_{nm}(N_{zsyr}^{\text{rep}})\mid  N_{zsyr}\},
\end{eqnarray*}
which can be approximated by the posterior draws from the Gibbs Sampler.

\section{\bf Evaluation of surrogate based on principal stratification}
\label{sec::evalution-principal-surrogacy}
In the previous section, we proposed approaches to identifying the PSACEs using multiple trials. 
In this section, our goal is to predict the  treatment  effect on the endpoint both quantitatively and qualitatively, based on the treatment effect  on the surrogate in a new trial without observing the endpoint. We consider two cases for applying the principal surrogate to a new trial in which the endpoint is not observed.
\begin{enumerate}
[{Case}~1.]
\item
The surrogate is applied to a new population
where the treatment is the same, but the distribution of principal strata differs from the validation trials. For example, in validation trials, we study the effect of drug A on cardiovascular disease, and use an indicator for whether the level of cholesterol increases as a candidate surrogate. Based on the PSACEs estimated  from validation trials, we use the surrogate to predict the effect of the drug A on the disease of patients in a different population.

\item
The  surrogate is applied to a new population
where both the treatment and the distribution of principal strata may be different from those in the validation trials. For example, based on the PSACEs for drug A  estimated  from validation trials,  we use the surrogate to predict the effect of  a new drug B on the disease  of patients in a different population.
\end{enumerate}

For quantitative evaluation, we assume that the PSACEs are the same across the validation trials and the new trial, then try to calculate the treatment effect on the endpoint.
For qualitative evaluation, we relax
this assumption so that the PSACEs have the same signs across the validation trials and the new trial. Then we try to get the sign of the  the treatment effect on the endpoint.

\citet{vanderweele2013surrogate} points out that a principal surrogate
satisfying the causal necessity ($ACE_{u}=0$ for $u=ss,\bar{s}\bar{s}$) and the causal sufficiency ($ACE_{u}\neq 0$ for $u=s\bar{s},\bar{s}s$) may not avoid the surrogate paradox.
In the following example, we further illustrate that
a principal surrogate verified in a validation trial may not be used to correctly evaluate
the treatment effect on the unobserved endpoint
in a new trial, even if the new trial and the validation trials
have the same PSACEs.

 \begin{example}
 \label{ex:2}
Consider a  surrogate $S$
which is evaluated by validation trial $r$.
The surrogate satisfies the causal necessity with
$ACE_{ss}=ACE_{\bar{s}\bar{s}}=0$,
and the causal sufficiency with
$ACE_{s\bar{s}}=0.4$ and $ACE_{\bar{s}s}=-0.6$.
Suppose that the distribution of principal strata
in the validation trial
is $\pi_{ss,r}=\pi_{\bar{s}\bar{s},r}=0.2$, $\pi_{s\bar{s},r}=0.4$ and $\pi_{\bar{s}s,r}=0.2$.
Then we have $E\{S(1)-S(0)\mid R=r \}=0.2>0$ and
$E\{Y(1)-Y(0)\mid R=r \} = 0.04>0$,
which means that the surrogate paradox is avoided in
validation trial $r$.
Now suppose that a new trial $R=r'$  has the same PSACEs as the validation trial.
But the distribution of principal strata differs: $\pi_{ss,r'}=0.1$, $\pi_{\bar{s}\bar{s},r'}=0.2$, $\pi_{s\bar{s},r'}=0.4$, and $\pi_{\bar{s}s,r'}=0.3$.
Then we have $E\{S(1)-S(0)\mid R=r' \}=0.1>0$ but
$E\{Y(1)-Y(0)\mid R=r' \}= -0.02<0$. 
The surrogate paradox arises,
and thus we cannot use the surrogate $S$
to correctly evaluate the treatment effect on the endpoint $Y$
in the new trial.
\end{example}

Now we discuss the conditions to avoid the surrogate
paradox when a validated surrogate
is applied to a new trial with a new treatment or a new population.
Let $ACE^S_r = E\{S(1) - S(0) \mid R=r\}$ and $ACE^Y_r = E\{Y(1) - Y(0) \mid R=r\}$.
Without loss of generality, we assume $ACE^S_r>0$. For the case with $ACE^S_r<0$, we can redefine the surrogate as $S^*=1-S$.
Below we give the relationships among the average causal effects 
$ACE^S_r, ACE^Y_r$ and $ACE_{ur}$ for trial $r$.

\begin{proposition}
\label{prop:sur:quanti}
For trial $r$, we assume that the surrogate satisfies the causal necessity.
\begin{enumerate}
\item[(i)]
With monotonicity, we have $ACE^Y_r = ACE^S_r \times ACE_{s\bar{s},r}$.

\item[(ii)] 
Without monotonicity, suppose $ACE^S_r>0$, we have the lower and upper bounds of $ACE^Y_r$: if $ACE_{s\bar{s},r}+ACE_{\bar{s}s,r} \geq 0$,
$$
ACE^S_r \times ACE_{s\bar{s},r} \leq ACE^Y_r \leq 
(ACE_{s\bar{s},r}+ACE_{\bar{s}s,r}) / 2 +  ACE^S_r\times  (ACE_{s\bar{s},r} - ACE_{\bar{s}s,r}) / 2,
$$
and otherwise
$$
(ACE_{s\bar{s},r}+ACE_{\bar{s}s,r}) / 2 +  ACE^S_r\times  (ACE_{s\bar{s},r} - ACE_{\bar{s}s,r}) / 2 \leq ACE^Y_r  \leq ACE^S_r \times ACE_{s\bar{s},r}.
$$

\end{enumerate}
\end{proposition}

Proposition \ref{prop:sur:quanti} shows the relationships between
the treatment effect on the surrogate and the treatment effect on the endpoint, which immediately give us 
the following implication relationships when the endpoint $Y$ is unobservable.

\begin{corollary}
 \label{coro:sur:quanli}
For trial $r$, we assume that the surrogate satisfies the causal necessity  and $ACE_{s\bar{s},r}> 0$.
\begin{enumerate}
\item[(i)] With monotonicity,
 $ACE^S_r> 0$     (or $=0$) implies $ACE^Y_r >0$ (or $=0$).

\item[(ii)] Without monotonicity,
$ACE_{s\bar{s},r} + ACE_{\bar{s}s,r} \geq 0$ and $ACE^S_r > 0$ imply $ACE^Y_r>0$.
\end{enumerate}
\end{corollary}

In  Corollary \ref{coro:sur:quanli}{\it (i)}, with monotonicity, a principal surrogate satisfying the causal sufficiency
can avoid the surrogate paradox.
Without monotonicity,
the above result {\it (i)} does not hold
since the causal effects on the endpoint $Y$
may be different in the principal strata $u=s\bar{s}$ and $\bar{s}s$.
Thus, we require the condition $ACE_{s\bar{s},r}+ACE_{\bar{s}s,r} \geq 0$ in {\it (ii)} of Corollary \ref{coro:sur:quanli}, i.e.,
the positive causal effect on the endpoint in
stratum $s\bar{s}$ can offset the negative causal effect
on the endpoint in stratum $\bar{s}s$. 

Proposision \ref{prop:sur:quanti} and Corollary \ref{coro:sur:quanli} can help us to assess the treatment effect on the endpoint  in a new trial both quantitatively and qualitatively. Below we illustrate this with Examples 3 and 4 for Cases 1 and 2 respectively.

\begin{example}\label{eg--3}
Suppose that in the validation trial, we studied drug A and  obtained estimates  $ACE_{ss,r}=ACE_{\bar{s}\bar{s},r}=0$, $ACE_{s\bar{s},r}=0.5$ and $ACE_{\bar{s}s,r}=-0.4$. For Case 1, we are aiming to evaluate the effect of the same drug on the disease  for a different population in a new trial $r'$. In trial $r'$, we obtain $E\{S(1) \mid R=r'\}=0.6 $ and $E\{S(0) \mid R=r'\}=0.4$. If we assume that the PSACEs in the new trial are the same as  those in the validation trial and monotonicity holds in the new trial, Proposition \ref{prop:sur:quanti}(i) implies that 
\begin{eqnarray*}
ACE^Y_{r'}= ACE^S_{r'} \times ACE_{s\bar{s},r'}
= ACE^S_{r'} \times ACE_{s\bar{s},r}=(0.6-0.4) \times 0.5=0.10 > 0.
\end{eqnarray*}
If monotonicity fails in the new trial, we can still calculate the bounds for $ACE^Y_{r'}$ according to Proposition \ref{prop:sur:quanti}(ii):
\begin{eqnarray}
 ACE^Y_{r'} &\leq& (ACE_{s\bar{s},r'}+ACE_{\bar{s}s,r'}) / 2 +  ACE^S_{r'}
\times  (ACE_{s\bar{s},r'} - ACE_{\bar{s}s,r'}) / 2  \nonumber \\
&=&(0.5-0.4)/2+(0.6-0.4)\times (0.5+0.4)/2 = 0.14, \label{eq::bound-1}\\
  ACE^Y_{r'} &\geq& ACE^S_{r'} \times ACE_{s\bar{s},r'} = (0.6-0.4) \times 0.5 =0.10.  \label{eq::bound-2}
\end{eqnarray}
If we only assume that $ACE_{ss,r'}=ACE_{\bar{s}\bar{s},r'}=0$ and $ACE_{s\bar{s},r'}+ACE_{\bar{s}s,r'}>0$ in the new trial as in the validation trial, then $ACE_{r'}^S=0.2>0$ and $ACE_{s\bar{s},r'}=0.5>0$ imply $ACE^Y_{r'}>0$, according to Corollary \ref{coro:sur:quanli} with or without monotonicity.
\end{example}

\begin{example}
Suppose that in the validation trial, we studied  drug A  and  estimated  $ACE_{ss,r}=ACE_{\bar{s}\bar{s},r}=0$, $ACE_{s\bar{s},r}=0.5$ and $ACE_{\bar{s}s,r}=-0.4$. For Case 2, we are aiming to evaluate the effect of another drug B  on the disease for the population in a new trial $r'$.  Unlike Case 1, we need to define the new drug as $Z'$, the principal stratification  as $U'=(S(Z'=1),S(Z'=0))$, and the PSACEs as $ACE_{ur'}=E\{Y(Z'=1)-Y(Z'=0)\mid U'=u, R=r'\}$.  In the new trial, we obtain $E\{S(1) \mid R=r'\}=0.6 $ and $E\{S(0) \mid R=r'\}=0.4$.  Similar to Example 3, if  we assume that the PSACEs in the new trial are the same as those in the validation trial and monotonicity holds in the new trial, we can calculate $ACE_{r'}^Y=0.1$. If monotonicity fails, we can obtain the bounds (\ref{eq::bound-1}) and (\ref{eq::bound-2}) for $ACE_{r'}^Y$.
If we assume $ACE_{ss,r'}=ACE_{\bar{s}\bar{s},r'}=0$ and $ACE_{s\bar{s},r'}+ACE_{\bar{s}s,r'}>0$ as in Example \ref{eg--3}, we can deduce that $ACE^Y_{r'}>0$ according to Corollary \ref{coro:sur:quanli}. 
\end{example}

Assuming the same values of PSACEs enables us to quantitatively evaluate the surrogate, and only assuming the same signs of PSACEs still allows us for qualitative evaluation. Note that we  give only the sufficient conditions to evaluate the surrogate in a new trial. 
The plausibilities of these assumptions depend on subject knowledge and experts' opinions.

\section{\bf Simulation studies}
\label{sec::simulation-studies}
In this section, we conduct simulation studies to evaluate finite sample performances of our proposed approaches, under both correctly specified and misspecified models.

\subsection{\bf Estimation under the models with homogeneity}

We show the simulation results of the MLEs and the credible intervals of the PSACEs with and without monotonicity.
The $95\%$ credible intervals are obtained by the Gibbs Sampler with $20,000$ iterations and the first $4,000$ iterations as the burn-in period.
We repeat $200$ times to get the coverage proportions of the true PSACEs for each setting.
We use three different sample sizes, and let $N$ denote the average sample size for each trial, i.e., $N= \text{total sample size}/N_R$.
In Figures \ref{fig:sim1:mon} and \ref{fig:sim1:non}, we show the biases and root mean squared errors (RMSEs) of the MLEs and the coverage proportions of the posterior credible intervals only for $ACE_{ss}$.
The results for other principal strata are similar and shown in the online Appendix.

First, with monotonicity,
we generate $R$ from a categorical distribution with $p_r =1/N_R$ for all $r$. We generate $Z$ from Bernoulli distributions conditional on $R$
with different conditional probabilities, $\bm{\alpha}$, allowing for different treatment assignment probabilities for different trials. We generate $U$ from categorical distributions conditional on $R$ with probabilities $(\pi_{ss,r}, \pi_{s\bar{s},r}, \pi_{\bar{s}\bar{s},r})$ for trial $r$. In order to satisfy Assumption \ref{assume:5}, we choose $U \mid R=r$ to depend on $R=r$.  We generate three scenarios with $N_R=2,3$ and $5$, and present the true values of $(\pi_{ss,r},\pi_{s\bar{s},r},\pi_{\bar{s}\bar{s},r})$ and $\bm{\alpha}$ in the upper panel of Table \ref{tab:sim}.
We generate $Y$ from Bernoulli distributions with
the following conditional probabilities given $Z$ and $U$:
$$
(\delta_{1, ss}, \delta_{0, ss}, \delta_{1, s\bar{s}}, \delta_{0, s\bar{s}}, \delta_{1, \bar{s}\bar{s}}, \delta_{0, \bar{s}\bar{s}}) = (0.8, 0.5, 0.7, 0.3, 0.6, 0.1),
$$
with PSACEs $ACE_{ss}=0.3$, $ACE_{s\bar{s}}=0.4$, and $ACE_{\bar{s}\bar{s}}=0.5$.

\begin{table}
\caption{\label{tab:sim}True parameters in simulations. We present the values of $\pi_{ur}=P(U=r \mid R=r)$ and $\alpha_r=P(Z=1 \mid R=r)$ in each scenario. The upper panel shows the true parameters with monotonicity, and
the lower panel shows the true parameters without monotonicity.}

\tabcolsep=0.11cm
\scalebox{0.9}{
\begin{tabular}{ccccccccc}\hline
\multicolumn{9}{c}{With monotonicity} \\ 
                 &\multicolumn{2}{c}{$N_R=2$} &&\multicolumn{2}{c}{$N_R=3$} &&\multicolumn{2}{c}{$N_R=5$} \\ 
                 \cline{2-3} \cline{5-6} \cline{8-9}
                 & $(\pi_{ss,r},\pi_{s\bar{s},r},\pi_{\bar{s}\bar{s},r})$ & $\alpha_r$ && $(\pi_{ss,r},\pi_{s\bar{s},r},\pi_{\bar{s}\bar{s},r})$ & $\alpha_r$ & & $(\pi_{ss,r},\pi_{s\bar{s},r},\pi_{\bar{s}\bar{s},r})$ & $\alpha_r$  \\
$r=1$        &  $(0.7,0.2,0.1)$  & 0.4  && $(0.8,0.1,0.1)$  &    0.4      & &  $(0.8,0.1,0.1)$  &  0.3\\
$r=2$        &  $(0.1,0.2,0.7)$  & 0.6 &&  $(0.1,0.8,0.1)$  &    0.5      & &  $(0.6,0.3,0.1)$  &0.4\\
$r=3$        &                           &       &&   $(0.1,0.1,0.8)$ &    0.6     &&   $(0.3,0.2,0.5)$ &0.5\\
$r=4$        &                           &       &&              &                            &&    $(0.1,0.3,0.6)$  &0.6\\
$r=5$        &                            &      &&               &                          &&  $(0.1,0.1,0.8)$     &0.7\\ \hline
\multicolumn{9}{c}{Without monotonicity} \\
                 &\multicolumn{2}{c}{$N_R=3$} &&\multicolumn{2}{c}{$N_R=4$} &&\multicolumn{2}{c}{$N_R=5$} \\ 
                  \cline{2-3} \cline{5-6} \cline{8-9}
                 & $(\pi_{ss,r},\pi_{s\bar{s},r},\pi_{\bar{s}\bar{s},r},\pi_{\bar{s}s,r})$ & $\alpha_r$ && $(\pi_{ss,r},\pi_{s\bar{s},r},\pi_{\bar{s}\bar{s},r},\pi_{\bar{s}s,r})$ & $\alpha_r$  && $(\pi_{ss,r},\pi_{s\bar{s},r},\pi_{\bar{s}\bar{s},r},\pi_{\bar{s}s,r})$ & $\alpha_r$  \\
$r=1$        &  $(0.6,0.2,0.1,0.1)$  & 0.4  && $(0.6,0.2,0.1,0.1)$  &    0.4       &&  $(0.6,0.2,0.1,0.1)$  &  0.3\\
$r=2$        &  $(0.1, 0.6,0.2,0.1)$  & 0.5 &&  $(0.1, 0.6,0.2,0.1)$  &    0.5       &&  $(0.1, 0.6,0.2,0.1)$  &0.4\\
$r=3$        &    $(0.1,0.1,0.6,0.2)$   &  0.6  &&  $(0.1,0.1,0.6,0.2)$ &    0.6     &&   $(0.3,0.2,0.3,0.2)$ &0.5\\
$r=4$        &        &   &&  $(0.2,0.3,0.2,0.3)$ &         0.7              &&    $(0.4,0.1,0.4,0.1)$ &0.6\\
$r=5$        &    &  & &                            &                                         & &  $(0.1,0.1,0.6,0.2)$     &0.7\\ \hline
\end{tabular}
}
\label{default}
\end{table}%

We show the simulation results in Figure \ref{fig:sim1:mon},
where the numbers ``2,'' ``3'' and ``5'' denote the numbers of trials.
The biases are small for all scenarios,
the RMSEs decrease as the sample size and the number of trials increase,
and the coverage proportions are close to 95\% for all scenarios. Comparing the 6th point with the  7th point in the second subplot of Figure \ref{fig:sim1:mon}, we can see that the RMSE of the case with $(N_R=2$, $N=500)$ is smaller than that of the case with $(N_R=5, N=200)$. The total sample sizes for these two cases are the same, while the case of the 6th point has more  parameters (more $\pi_{ur}$'s with more trials). The simulation is consistent with our intuition that more  parameters will result in a larger RMSE.

Next, without monotonicity, we generate $R$ from a categorical distribution with $p_r =1/N_R$ for all $r$. We generate $Z$ from Bernoulli distributions with conditional probabilities $\bm{\alpha}$, and $U$ from categorical distributions both conditional on $R$ with probabilities $(\pi_{ss,r}, \pi_{s\bar{s},r}, \pi_{\bar{s}\bar{s},r},\pi_{\bar{s}s,r})$ for trial $r$. For this case without  monotonicity,
the necessary condition for identifiability
is $N_R\geq 3$.
Thus we generate three scenarios with $N_R=3,4$ and $5$, and present the true values of $(\pi_{ss,r},\pi_{s\bar{s},r},\pi_{\bar{s}\bar{s},r},\pi_{\bar{s}s,r})$ and $\bm{\alpha}$ in the lower panel of Table \ref{tab:sim}.
We generate $Y$ from Bernoulli distributions with conditional probabilities given $U$ and $Z$ as
$$
(\delta_{1, ss}, \delta_{0, ss}, \delta_{1, s\bar{s}}, \delta_{0, s\bar{s}}, \delta_{1, \bar{s}\bar{s}}, \delta_{0, \bar{s}\bar{s}}, \delta_{1, \bar{s}s}, \delta_{0, \bar{s}s})
= (0.8, 0.5, 0.7, 0.3, 0.6, 0.1, 0.5, 0.2),
$$
with PSACEs $ACE_{ss}=0.3$, $ACE_{s\bar{s}}=0.4$, $ACE_{\bar{s}\bar{s}}=0.5$, and $ACE_{\bar{s}s}=0.3$.

We show the simulation results in Figure \ref{fig:sim1:non},
where ``3,'' ``4'' and ``5'' denote the numbers of trials.
The biases and the coverage proportions are similar to
the cases with monotonicity.
But the RMSEs become a little larger than those with 
monotonicity.
The performances of our approaches are quite promising for the sample size
$N=500 \times N_R$, comparable to the ACCTs data set in Section \ref{sec::application-ACCT}.

\subsection{\bf Simulations under the models without homogeneity}

We conduct simulation studies
when the homogeneity assumption is false with $\delta_{zur}$ dependent on $r$.
Let $d$ be a measure of heterogeneity.
We simulate three cases for different levels of heterogeneity with $d= 0.01, 0.025$ and $0.05$.
For each case, we conduct simulation studies when $N_R=3$ both with and without monotonicity.

With monotonicity, we first generate $(Z,U)$ the same as 
the scenario with $N_R=3$ in the upper panel of Table \ref{tab:sim}.
To allow heterogeneity, we let $\delta_{zu1}=\mu_{zu}-(-1)^{z}d, \delta_{zu2}=\mu_{zu}$ and $ \delta_{zu3}=\mu_{zu}+(-1)^{z}d$, with $\delta_{zur}$ varying with $r$.
The corresponding mean vector is
$$
(\mu_{1, ss}, \mu_{0, ss},  \mu_{1, s\bar{s}}, \mu_{0, s\bar{s}}, \mu_{1, \bar{s}\bar{s}}, \mu_{0, \bar{s}\bar{s}} ) = ( 0.8, 0.5, 0.7, 0.3, 0.6, 0.1 ).
$$
Then we have $ACE_{u1} = \mu_{1u}-\mu_{0u}+2d$, $ACE_{u2} = \mu_{1u}-\mu_{0u}$ and $ACE_{u3} = \mu_{1u}-\mu_{0u}-2d$. A larger $d$ results in a larger violation of the homogeneity assumption.

Similarly, without monotonicity,
we first generate $(Z,U)$ the same as the scenario with $N_R=3$ in the lower panel of Table \ref{tab:sim}.
Then
we let $\delta_{zu1}=\mu_{zu}-(-1)^{z}d, \delta_{zu2}=\mu_{zu}$ and $\delta_{zu3}=\mu_{zu}+(-1)^{z}d$, with the corresponding mean vector
$$
( \mu_{1, ss}, \mu_{0, ss},  \mu_{1, s\bar{s}}, \mu_{0, s\bar{s}}, \mu_{1, \bar{s}\bar{s}}, \mu_{0, \bar{s}\bar{s}}, \mu_{1, \bar{s}s}, \mu_{0, \bar{s}s} )
 = (0.8, 0.5, 0.7, 0.3, 0.6, 0.1, 0.5, 0.2 ).
$$
Then we have $ACE_{u1} = \mu_{1u}-\mu_{0u}+2d$, $ACE_{u2} = \mu_{1u}-\mu_{0u}$, and $ACE_{u3} = \mu_{1u}-\mu_{0u}-2d$.

\begin{figure}
  \centering
   \subfigure[Correctly specified model  with monotonicity]{
    \label{fig:sim1:mon} 
\scalebox{1}[0.95]{\includegraphics[width=\textwidth]{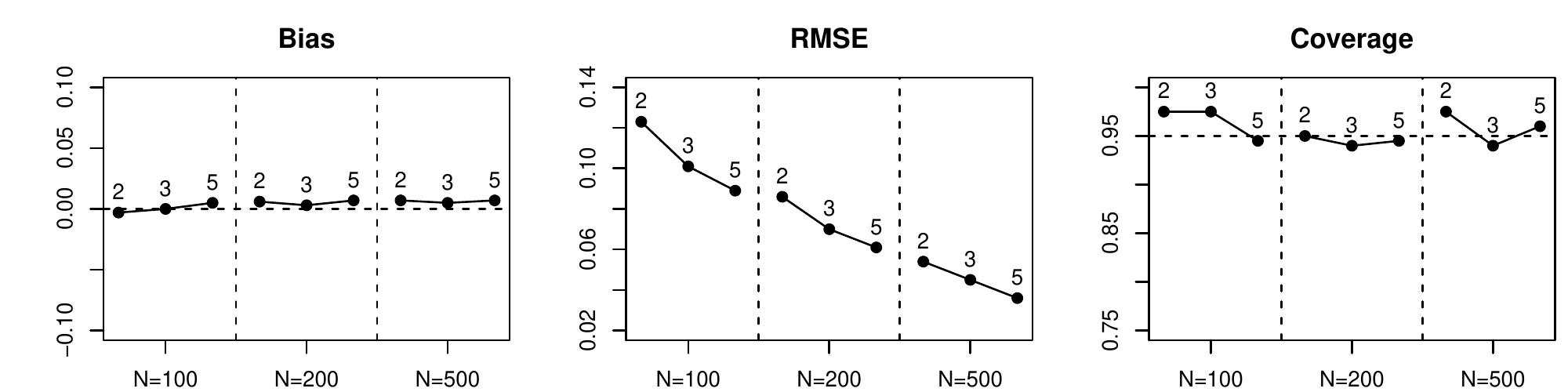}}} \subfigure[Correctly specified model  without monotonicity]{
    \label{fig:sim1:non} 
\scalebox{1}[0.95]{\includegraphics[width=\textwidth]{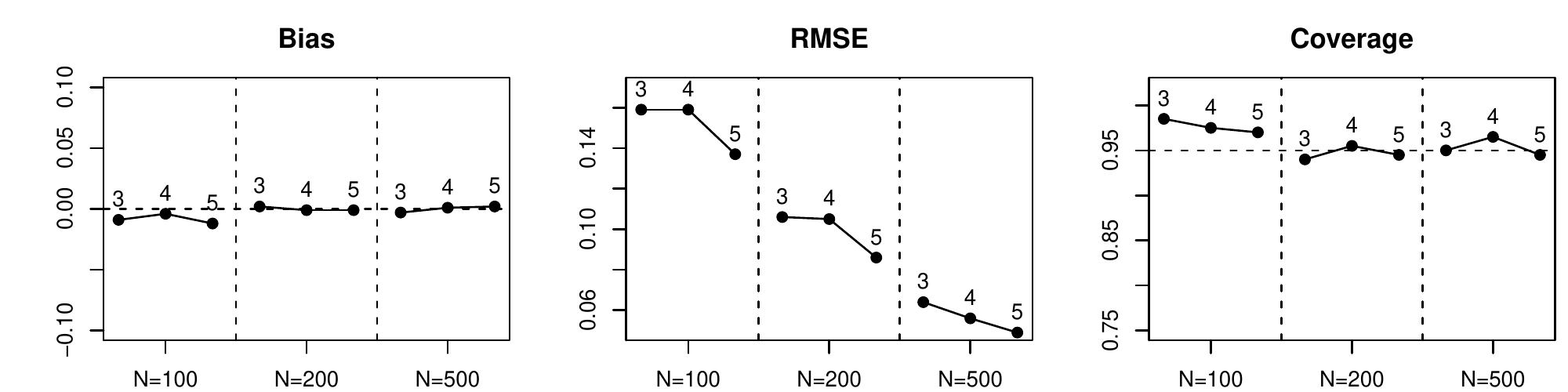}}}
  \subfigure[Misspecified model with monotonicity]{
    \label{fig:sim2:mon} 
\scalebox{1}[0.95]{\includegraphics[width=\textwidth]{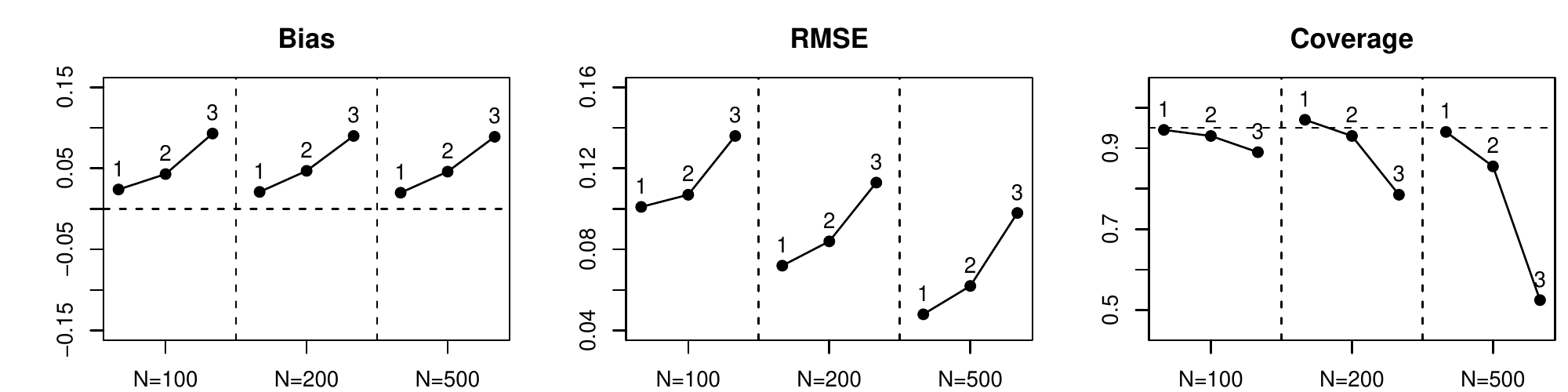}}} \subfigure[Misspecified model without monotonicity]{
    \label{fig:sim2:non} 
\scalebox{1}[0.95]{\includegraphics[width=\textwidth]{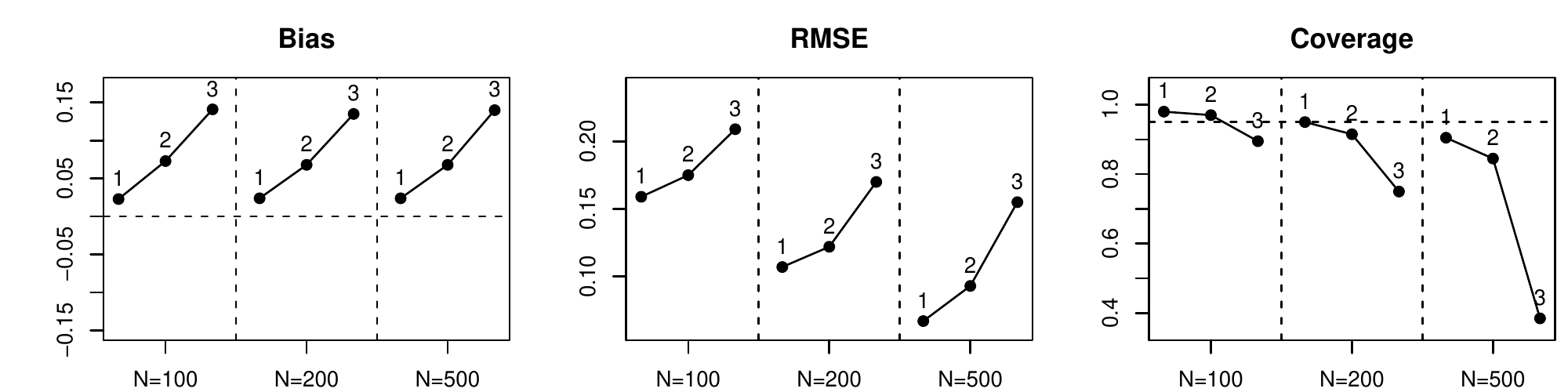}}}
  \caption{Simulations for homogeneity (a) and (b) and heterogeneity (c) and (d). Each subgraph presents the bias, RMSE, coverage proportions of the $95\%$ credible intervals of $ACE_{ss}$. Nine combinations of $N_R$ and $N$ are shown in (a) and (b); and nine combinations of $N$ and $d$ (``1'' for $.01$, ``2'' for $.025$, ``3'' for $.05$) are shown in (c) and (d).}
  \label{fig:sim2} 
\end{figure}

For all the evaluations, we use the mean parameters $\mu_{zu}$ as the ``true parameters.''
Figures \ref{fig:sim2:mon} and \ref{fig:sim2:non}
show the results for $ACE_{ss}$ including the biases, the RMSEs of the MLEs and the coverage proportions of the posterior credible intervals,
where  ``1,'' ``2'' and ``3'' correspond to the cases of
$d=0.01$, $0.025$ and $0.05$, respectively.

Comparing Figures \ref{fig:sim1:mon} and \ref{fig:sim1:non} under the homogeneity assumption   with
Figures \ref{fig:sim2:mon} and \ref{fig:sim2:non} without the homogeneity assumption, we find that the biases increase as $d$ increases, and do not decrease as the sample size grows. 
The RMSEs increase and the coverage proportions of $95\%$ credible intervals decrease as $d$ increases, especially for large sample sizes.
This means that the point and interval estimates are sensitive to the heterogeneity of trials. Thus, we will propose an approach for sensitivity analysis, and apply it to our real application in Section 7.

\section{\bf Application to the ACCTs data}
\label{sec::application-ACCT}
For the ACCTs, the OS with five-year follow-up (survival status within 5 years) is used as
the endpoint in most earlier papers to evaluate a particular treatment regimen. However, this endpoint requires five-year follow-up. The goal of the ACCTs study is to explore whether DFS with three-year follow-up
(cancer reoccurred before 3 years) is a valid surrogate for the OS with five-year follow-up.
The data contain more than $20,000$ patients and $18$ randomized III clinical trials.  The period of trial enrollment spans from 1977 to 1999.
  The data of $10$ of the randomized trials are available from \citet{baker2012predicting}.  \citet{baker2012predicting} transformed the survival data into counts for binary outcomes. As a result, in each trial, we have  a contingency table of observed frequencies with three variables: the treatment ($Z$), DFS ($S$) and OS ($Y$). To illustrate our approaches, we assume all the trials have  the same treatment and control, although we do not know the exact treatment and control in each trial. This may influence the plausibility of the homogeneity assumption, and thus we will conduct a sensitivity analysis in Section 7.
We show all the bound analysis in the online Appendix, and find that the bounds are barely informative.
In this section, we apply the proposed approaches to the ACCTs data to obtain point identifications of the treatment effects on the endpoint OS within principal strata defined by the potential DFS. 
We postpone the sensitivity analysis to Section \ref{sec::sensitivity}.

\subsection{\bf Results with monotonicity}
We use the Gibbs Sampler to simulate the posterior distributions of the PSACEs, because it is very direct to obtain posterior credible intervals for the PSACEs from draws of the Gibbs Sampler. After 100,000 iterations with 50,000 used as burn-in, the Markov chains converge very well with Gelman--Rubin diagnostic statistics approximately equal to $1$ from five independent chains.
Figure \ref{fig:hist_mon} shows the histograms of the posterior draws of the PSACEs under monotonicity.
The posterior median of $ACE_{\bar{s}\bar{s}}$, $-0.008$, is small, and its $95\%$ credible interval covers zero. Although the posterior median of $ACE_{ss}$, $0.032$, is also small, its $95\%$ credible interval does not cover zero. If we would like to believe monotonicity, there exists a ``principal stratification direct effect'' of the treatment on the endpoint for stratum $ss$ but not for stratum $\bar{s}\bar{s}$, and thus
the causal necessity is not well satisfied.
Because $ACE_{s\bar{s}}$ is significantly nonzero with posterior median $0.583$ and $95\%$ credible interval $(0.235,0.787)$,  and the treatment is effective on the surrogate with posterior medians of $ACE_r^S$ greater than $0.022$ in all trials, we can conclude that it is also effective on the endpoint according to Corollary \ref{coro:sur:quanli} {\it (i)}.

\subsection{\bf Results without monotonicity}
\label{sec::results-without-mono}

\begin{figure}[htpb]
\centering
\subfigure[Posterior distributions of the PSACEs with monotonicity.]{
\label{fig:hist_mon}
\includegraphics[width=\textwidth]{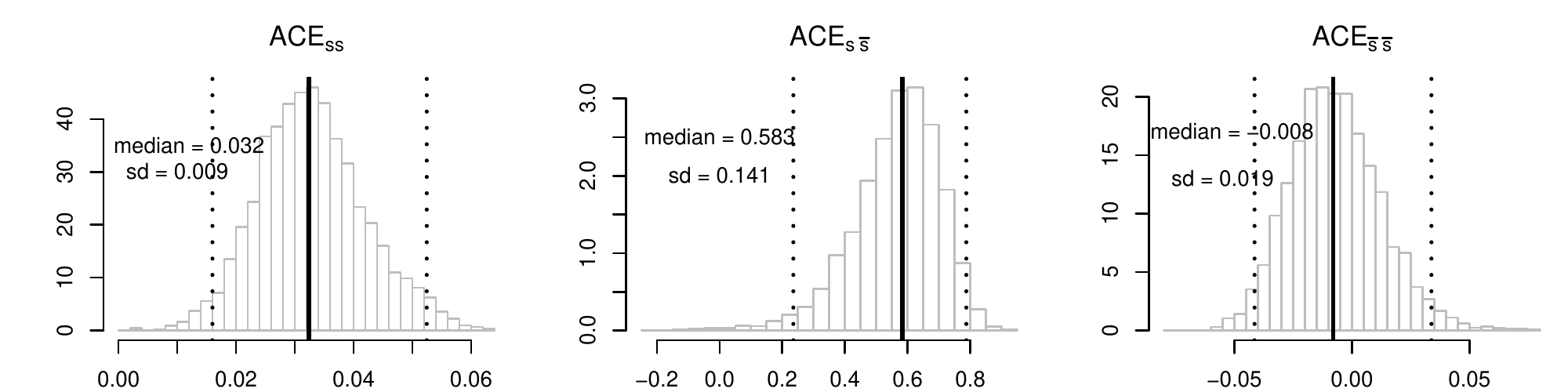}}
\subfigure[Posterior distributions of the PSACEs without monotonicity.]{
\label{fig:hist}
\includegraphics[width=\textwidth]{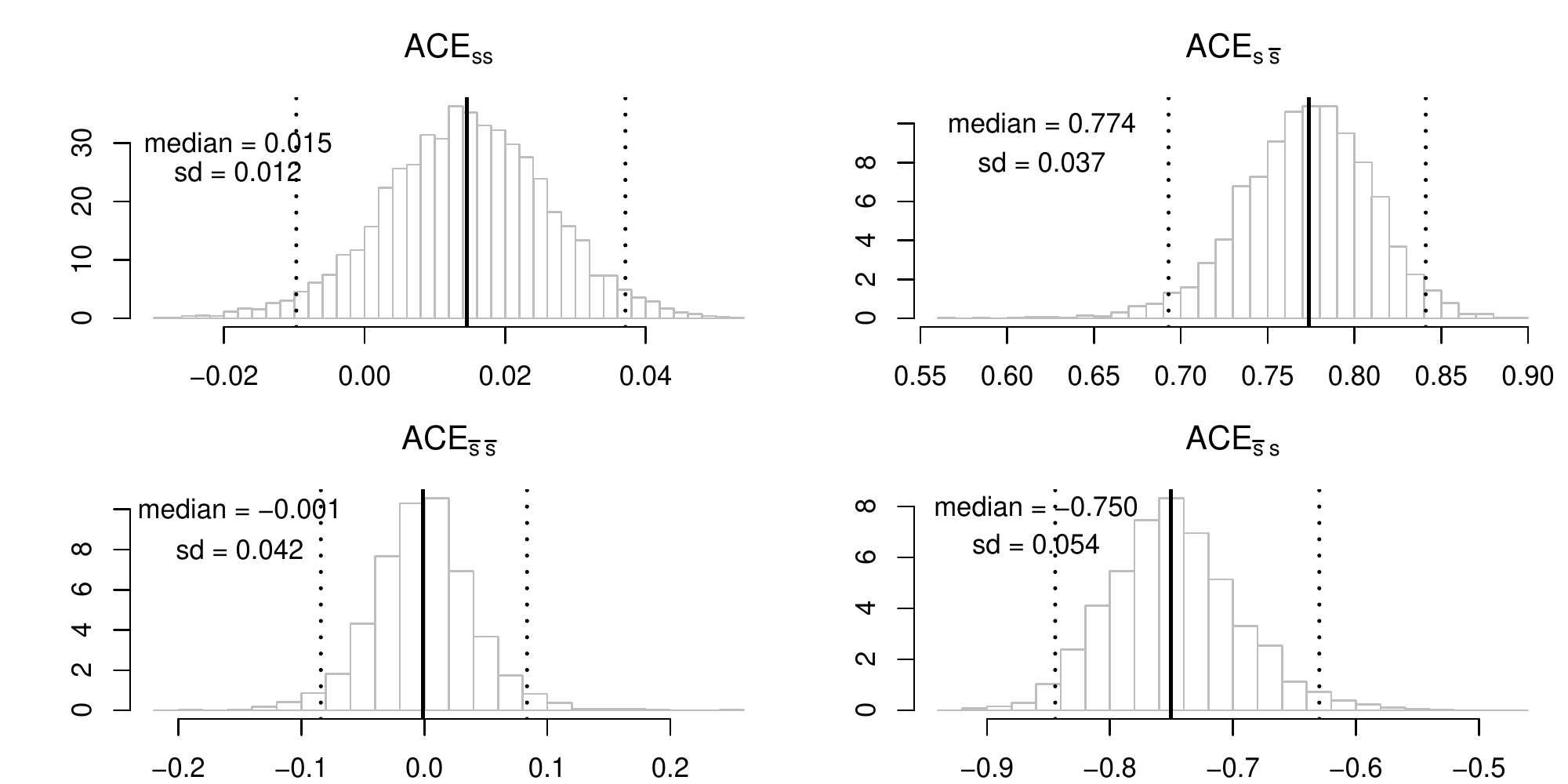}
}
\subfigure[Posterior distributions of the principal strata without monotonicity.]{
\label{fig:lh}
\includegraphics[width=1\textwidth]{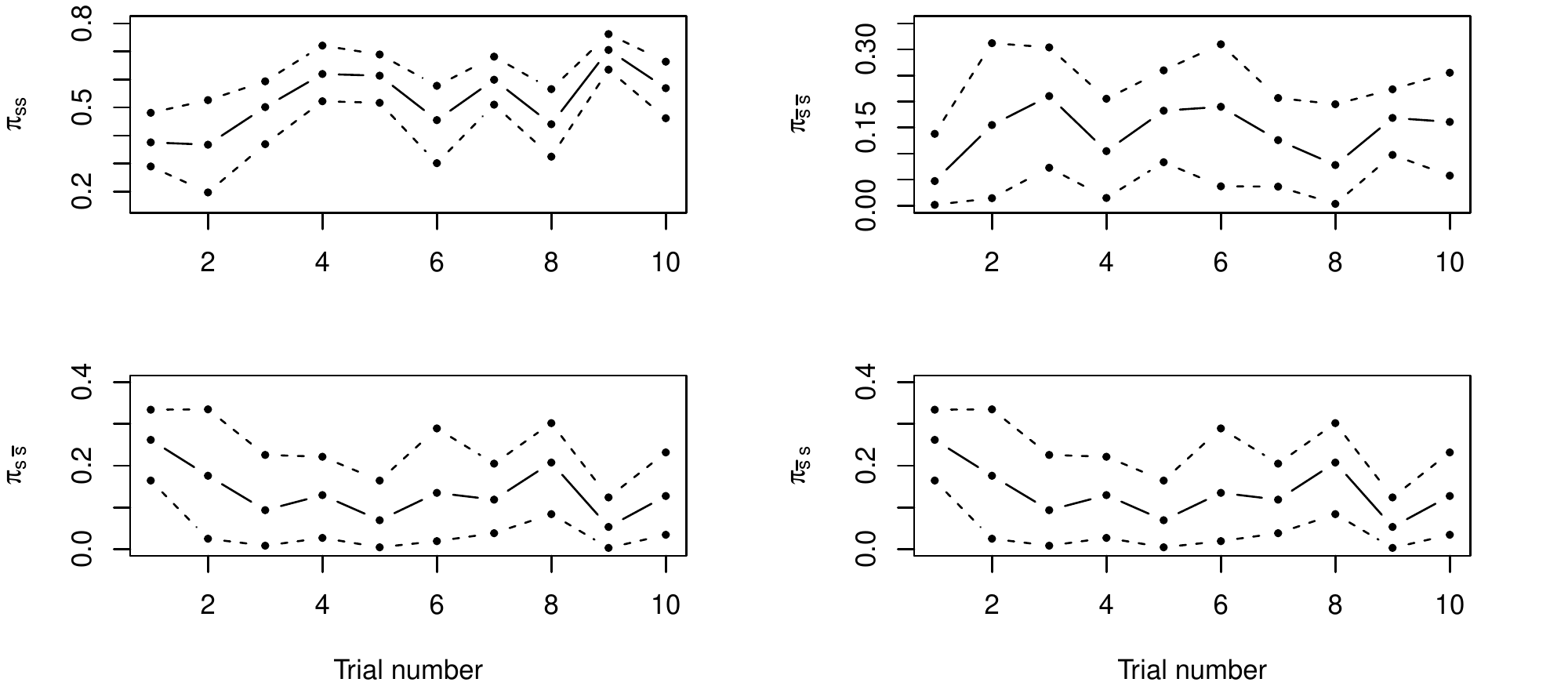}
}
\caption{The solid lines are the posterior medians, and the dotted lines are the posterior $2.5\%$ and $97.5\%$ quantiles.}
\end{figure}

Figure \ref{fig:hist} shows the histograms of posterior draws of the PSACEs without  monotonicity. First,
$ACE_{ss}$ and $ACE_{\bar{s}\bar{s}}$ are very close to zero with posterior medians $0.014$ and $0.001$ and $95\%$ credible intervals covering zero.
When the treatment does not affect the DFS with three-year follow-up, it will not affect the OS with five-year follow-up.
Second, $ACE_{s\bar{s}}$ and $ACE_{\bar{s}s}$ are significantly different from zero with posterior medians $0.774$ and $-0.750$ and $95\%$ credible intervals excluding zero.
Therefore, when the treatment affects the DFS with three-year follow-up, it will also affect the OS with five-year follow-up.
The results without monotonicity differ from the results with monotonicity in the posterior distribution of $ACE_{ss}$, which give us different interpretations about the causal mechanism.
From a practical point of view, we need to verify which model is more plausible based on the observed data.

\subsection{\bf Model comparison and checking}
We use the methods described in Section \ref{sec::computation} to perform model comparison.
The ppp's are $0.039$ and $0.470$ under the models with and without monotonicity, respectively.
Analogously, the $p$-values from the likelihood ratio tests are $0.015$ and $ 0.140 $ under the models with and without monotonicity, respectively. Therefore, the model with monotonicity is not compatible with the observed data, and it is rejected by both frequentists' and Bayesian methods. Furthermore, the lower right subplot of Figure \ref{fig:lh} shows that the posterior medians and $95\%$ credible intervals of $\pi_{\bar{s}s,r}$ under the model without monotonicity are far from zero in many trials.
On the contrary, the large $p$-values from both the likelihood ratio test and the posterior predictive check indicate very good fit of the model without monotonicity to the observed data.
Although previous methods often assumed monotonicity, our tests reject it. Therefore, our methods can be used in a broader scope of applications and are more credible when monotonicity fails.

\subsection{\bf Evaluation of principal surrogate in the ACCTs data}

According to the discussion  in Section 5 about model checking, we believe that the model without monotonicity is more credible, and thus all our discussion below is based on the results in Section \ref{sec::results-without-mono} from the model without monotonicity.

In Section \ref{sec::results-without-mono}, we obtain that the posterior medians of $ACE_{ss}$ and $ACE_{\bar{s}\bar{s}}$ are very close to zero, and both of their $95\%$ credible intervals cover zero, which means that the candidate surrogate satisfies the causal necessity very well.  Therefore, our analysis verifies the causal necessity assumption in Baker et al. (2012).
In addition, the posterior medians of $ ACE_{s\bar{s}}$ and $ ACE_{\bar{s}s}$ are  $0.772$ and  $- 0.748$, respectively, with both of their $95\%$ credible intervals excluding zero, showing that the surrogate also satisfies the causal sufficiency very well.  We then have approximately $ ACE_{ss}  = 0, ACE_{\bar{s}\bar{s}} =  0, ACE_{s\bar{s}} > 0$, and $ACE_{s\bar{s}}+ACE_{\bar{s}s}  \geq 0$. Thus, according to Corollary 1$(ii)$, the surrogate paradox can be avoided by  using this surrogate, i.e., if the treatment has a positive causal effect on the surrogate, it must have a positive causal effect on the endpoint. In fact, the treatment has a positive average causal effect on the surrogate (posterior medians of $ACE_r^S$ are greater than $0.015$ in all trials), and we can conclude that it has a positive average causal effect on the endpoint.

\section{\bf Sensitivity analysis without homogeneity}
\label{sec::sensitivity}

The homogeneity assumption  is crucial for the identifiability of PSACEs. It may be violated if different trials with different environments may affect the endpoint. Instead of assuming that $\delta_{zur}$ does not depend on $r$, we propose a Bayesian hierarchical model to account for the heterogeneity among $\delta_{zur}$ for different trials.
We keep the conditional distributions of $P(R)$, $P(Z \mid R)$ and $P(U \mid R)$ unchanged, but
assume the following hierarchical model for $P(Y \mid Z,U,R)$:
\begin{eqnarray*}
Y \mid  Z=z, U=u, R=r &\sim& \text{Bernoulli}(\delta_{zur}),\\
\text{logit}(\delta_{zur}) &\sim & N(\mu_{zu}, \sigma^2).
\end{eqnarray*}
In this model, deviation from the homogeneity assumption is characterized by the sensitivity parameter $\sigma$. When $\sigma = 0$, $\delta_{zur} =\delta_{zu}$ and thus the homogeneity assumption holds. When $\sigma > 0$, the homogeneity assumption is violated.
For example, when $\mu_{zu}=\text{logit}(0.3)$ and $\sigma=0.5$, the parameter $\delta_{zur}$ falls in the interval $(0.139, 0.533)$ with probability $0.95$; when $\mu_{zu}=\text{logit}(0.5)$ and $\sigma = 0.5$, the parameter $\delta_{zur}$ falls in the interval $(0.273, 0.727)$ with probability $0.95$.
Since the  parameter $\delta_{zur}$ is within the interval $[0, 1]$,
the above intervals imply quite large deviations away from the homogeneity assumption.
Therefore, we choose $0.5$ as the maximum value of the sensitivity parameter $\sigma$, and choose $0.05$ and $ 0.2$ as two moderate values of $\sigma$.
In our Bayesian analysis, we
choose the following priors:  $\{P(Z=1\mid R=r): r=1,\ldots,N_R\} \sim \text{Dirichlet}(1,\cdots,1)$, $ P(Z=1\mid R=r)$ $\sim U(0,1)$, $\{\pi_{ss,r}, \pi_{s\bar{s},r}, \pi_{\bar{s}\bar{s},r}, \pi_{\bar{s}s,r} \} \sim  \text{Dirichlet}(1,\cdots,1)$,  and $\mu_{zu} \sim U(-5,5)$. 
We use the prior of $\mu_{zu}$ for numerical stability, and also $[\text{logit}^{(-1)}(-5),\text{logit}^{(-1)}(5)]=[0.007,0.993] \approx [0,1]$ is not very restrictive.
The details of the Gibbs Sampler for the Bayesian hierarchical model are given in the online Appendix.

\begin{figure}[htpb]
\scalebox{1}[1.2]{\includegraphics[width=1\textwidth]{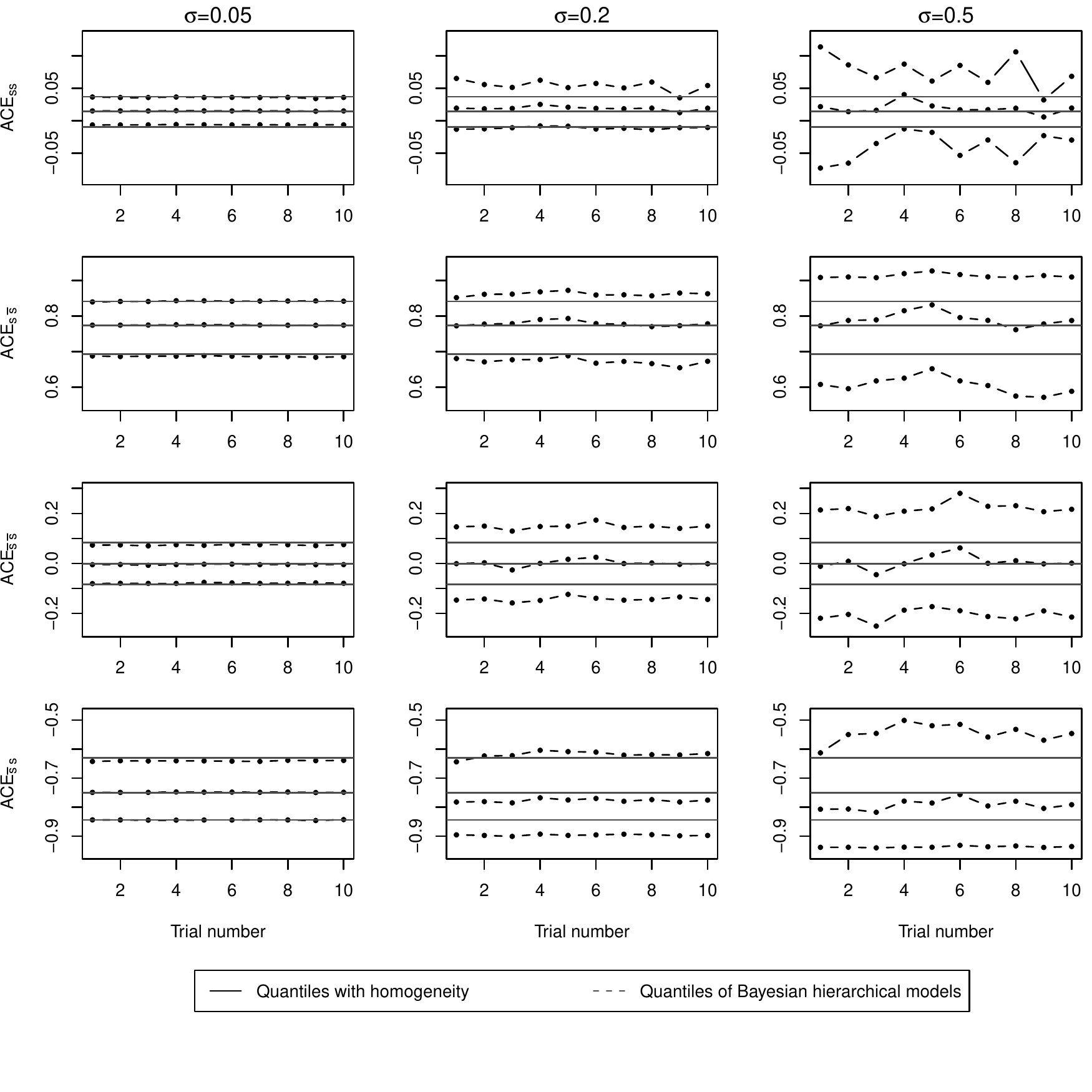}}
\caption{  
Sensitivity analysis for the ACCTs data.
The rows correspond to the four principal strata
and the columns correspond to the three different values of $\sigma$.
The dashed lines are the posterior $2.5\%$, $50\%$ and $97.5\%$ quantiles from Bayesian hierarchical models, and the solid grey lines are the posterior quantiles from the homogeneity model.}
\label{fig:sens}
\end{figure}

We reanalyze the ACCTs data using the Bayesian hierarchical models with results shown in Figure \ref{fig:sens}.
Four rows correspond to the four principal strata ($U=ss,s\bar{s}, \bar{s}s, \bar{s}\bar{s}$),
and three columns correspond to the three values of $\sigma$ ($\sigma=0.05,0.2,0.5$).
For each subgraph, we draw the $2.5\%$, $50\%$ and $97.5\%$ posterior quantiles of PSACEs obtained from both the homogeneity model and Bayesian hierarchical models. The posterior distributions of PSACEs are robust for these different values of the sensitivity parameter $\sigma$, although the intervals become wider as $\sigma$ increases.
The signs of $ACE_{s\bar{s}}$ and $ACE_{\bar{s}s}$
remain positive and negative respectively,
and the intervals of $ACE_{ss}$ and $ACE_{\bar{s}\bar{s}}$ still
 cover zero.
The results obtained by Bayesian hierarchical models
remain unchanged by removing the homogeneity assumption, which show the robustness of our results in Section \ref{sec::application-ACCT}.

\section{\bf Discussion}
\label{sec::discussion}

According to \citet{frangakis2002principal}, evaluation of surrogates requires identification of the causal effects within principal strata.
Because the potential outcomes of the surrogate endpoints under  treatment and control cannot be jointly observed, the value of the principal stratification variable is missing.
The monotonicity assumption excludes all the individuals with surrogate negatively affected by the treatment, which may not be biologically reasonable. The exclusion restriction assumption cannot be made when evaluating surrogates, because it implies the causal necessity and the validity of the causal necessity is a scientific question of interest.
Unfortunately, the identifiability of the PSACEs are jeopardized without monotonicity or exclusion restriction.
Although Bayesian analysis for weakly identified models still yields proper posterior distributions under proper priors, it may result in sensitive answers to the prior specification.

The ACCTs data contain multiple independent trials for evaluating the same surrogate for the same endpoint.
If we could make the homogeneity assumption across multiple trials,
we can remove the monotonicity and exclusion restriction assumptions but still guarantee identifiability of the PSACEs.
We find in the ACCTs data that
 both the ``causal necessity'' and the ``causal sufficiency'' hold, which imply that the DFS with three-year follow-up is a valid surrogate for OS with five-year follow-up.
We investigate the applicability of the principal surrogate in new trials on new populations or with new treatments, and show that the surrogate in the ACCTs data could be a useful surrogate under some conditions.
To remove the homogeneity assumption,
we further propose an approach based on Bayesian hierarchical models, and investigate the sensitivity to the deviations from the homogeneity assumption.
Within a reasonable range of the sensitivity parameter, the conclusions for the ACCTs data remain stable.

The framework we proposed could be applied to various settings involving post-treatment intermediate variables where the monotonicity and exclusion restriction assumptions are questionable. For example, in randomized experiments with non-compliance, the assignment of the treatment may influence the outcome directly. 
 In these cases, we should concern the plausibility of exclusion restriction assumptions versus the plausibility of homogeneity assumptions. Models with the homogeneity assumption may sometimes be more flexible in practice.  
 With the homogeneity assumption, we can test the monotonicity assumption. 
 Furthermore, even if the homogeneity assumption does not hold, we can still conduct a sensitivity analysis as in Section \ref{sec::sensitivity}.

\citet{follmann2006augmented}, \citet{qin2008assessing} and \citet{mattei2011augmented} proposed augmented designs for single trials to identify principal stratification causal effects using additional variables.
From the perspective of experimental design,
our approach can be viewed as an extension of \citet{follmann2006augmented}, which imposes a stronger monotonicity assumption for the treatment effect on the intermediate variable. 
 In our approaches, the number of trials $R$ can be more general, as long as it satisfies the assumptions for identification.  We can design experiments to create a variable $R$, satisfying these assumptions in order to identify the PSACEs.   For example, in an encouragement experiment that encourages the patients to take their assignments, we can use different types of encouragement for different groups of people.  Thus different groups will have different compliance behaviors, and the indicator for encouragement types then acts like $R$ in our paper.

Several generalizations are possible. 
First, it is relatively easy to deal with a continuous endpoint because we can identify the distributional causal effect by dichotomizing the endpoint \citep{ding2011identifiability}. However, it is non-trivial to deal with a continuous surrogate endpoint since the causal necessity and sufficiency cannot be defined by dichotomizing it. \citet{schwartz2011bayesian} set the stage for analyzing principal strata effects with continuous intermediate variables using a semiparametric Bayesian approach.
Second, we discuss the methods and present an application without missing data, and dealing with the missing data problem will be of interest in many real applications. Finally, generalizing our approaches to longitudinal data and time-to-event data is also of theoretical and practical interest.

\section*{Acknwoledgements}
We would like to thank the joint editor, the associate editor and
two referees for their valuable comments and suggestions which
greatly improved the previous version of this paper. This research
was supported by NSFC (11171365, 11331011) and 863 Program of China (2015AA020507).
Professor Carl N. Morris' pioneer work on hierarchical models motivated the second author to develop the method in Section \ref{sec::sensitivity}, and this paper is dedicated to him on the occasion of his retirement from Harvard University.

\bibliographystyle{chicago}
\footnotesize
\bibliography{PSACEsurrogate}

\begin{thebibliography}{}

\bibitem[\protect\citeauthoryear{Andrews}{Andrews}{2000}]{andrews2000inconsist%
ency}
Andrews, D.~W. (2000).
\newblock Inconsistency of the bootstrap when a parameter is on the boundary of
  the parameter space.
\newblock {\em Econometrica\/}~{\em 68}, 399--405.

\bibitem[\protect\citeauthoryear{Andrews and Guggenberger}{Andrews and
  Guggenberger}{2009}]{andrews2009validity}
Andrews, D.~W. and P.~Guggenberger (2009).
\newblock Validity of subsampling and plug-in asymptotic inference for
  parameters defined by moment inequalities.
\newblock {\em Econometric Theory\/}~{\em 25}, 669--709.

\bibitem[\protect\citeauthoryear{Angrist, Imbens, and Rubin}{Angrist
  et~al.}{1996}]{angrist1996identification}
Angrist, J.~D., G.~W. Imbens, and D.~B. Rubin (1996).
\newblock Identification of causal effects using instrumental variables (with
  discussion).
\newblock {\em J. Am. Statist. Ass.\/}~{\em 91}, 444--455.

\bibitem[\protect\citeauthoryear{Baker, Sargent, Buyse, and Burzykowski}{Baker
  et~al.}{2012}]{baker2012predicting}
Baker, S.~G., D.~J. Sargent, M.~Buyse, and T.~Burzykowski (2012).
\newblock Predicting treatment effect from surrogate endpoints and historical
  trials: an extrapolation involving probabilities of a binary outcome or
  survival to a specific time.
\newblock {\em Biometrics\/}~{\em 68}, 248--257.

\bibitem[\protect\citeauthoryear{Chen, Geng, and Jia}{Chen
  et~al.}{2007}]{chen2007criteria}
Chen, H., Z.~Geng, and J.~Jia (2007).
\newblock Criteria for surrogate end points.
\newblock {\em J. R. Statist. Soc. B\/}~{\em 69}, 919--932.

\bibitem[\protect\citeauthoryear{Cheng and Small}{Cheng and
  Small}{2006}]{cheng2006bounds}
Cheng, J. and D.~S. Small (2006).
\newblock Bounds on causal effects in three-arm trials with non-compliance.
\newblock {\em J. R. Statist. Soc. B\/}~{\em 68}, 815--836.

\bibitem[\protect\citeauthoryear{Chernozhukov, Lee, and Rosen}{Chernozhukov
  et~al.}{2013}]{chernozhukov2013intersection}
Chernozhukov, V., S.~Lee, and A.~M. Rosen (2013).
\newblock Intersection bounds: estimation and inference.
\newblock {\em Econometrica\/}~{\em 81}, 667--737.

\bibitem[\protect\citeauthoryear{Daniels and Hughes}{Daniels and
  Hughes}{1997}]{daniels1997meta}
Daniels, M.~J. and M.~D. Hughes (1997).
\newblock Meta-analysis for the evaluation of potential surrogate markers.
\newblock {\em Stat. Med.\/}~{\em 16}, 1965--1982.

\bibitem[\protect\citeauthoryear{Ding, Geng, Yan, and Zhou}{Ding
  et~al.}{2011}]{ding2011identifiability}
Ding, P., Z.~Geng, W.~Yan, and X.~H. Zhou (2011).
\newblock Identifiability and estimation of causal effects by principal
  stratification with outcomes truncated by death.
\newblock {\em J. Am. Statist. Ass.\/}~{\em 106}, 1578--1591.

\bibitem[\protect\citeauthoryear{Follmann}{Follmann}{2006}]{follmann2006augmen%
ted}
Follmann, D. (2006).
\newblock Augmented designs to assess immune response in vaccine trials.
\newblock {\em Biometrics\/}~{\em 62}, 1161--1169.

\bibitem[\protect\citeauthoryear{Frangakis and Rubin}{Frangakis and
  Rubin}{1999}]{frangakis1999addressing}
Frangakis, C.~E. and D.~B. Rubin (1999).
\newblock Addressing complications of intention-to-treat analysis in the
  combined presence of all-or-none treatment-noncompliance and subsequent
  missing outcomes.
\newblock {\em Biometrika\/}~{\em 86}, 365--379.

\bibitem[\protect\citeauthoryear{Frangakis and Rubin}{Frangakis and
  Rubin}{2002}]{frangakis2002principal}
Frangakis, C.~E. and D.~B. Rubin (2002).
\newblock Principal stratification in causal inference.
\newblock {\em Biometrics\/}~{\em 58}, 21--29.

\bibitem[\protect\citeauthoryear{Frumento, Mealli, Pacini, and Rubin}{Frumento
  et~al.}{2012}]{frumento2012evaluating}
Frumento, P., F.~Mealli, B.~Pacini, and D.~B. Rubin (2012).
\newblock Evaluating the effect of training on wages in the presence of
  noncompliance, nonemployment, and missing outcome data.
\newblock {\em J. Am. Statist. Ass.\/}~{\em 107}, 450--466.

\bibitem[\protect\citeauthoryear{Gallop, Small, Lin, Elliott, Joffe, and
  Ten~Have}{Gallop et~al.}{2009}]{gallop2009mediation}
Gallop, R., D.~S. Small, J.~Y. Lin, M.~R. Elliott, M.~Joffe, and T.~R. Ten~Have
  (2009).
\newblock Mediation analysis with principal stratification.
\newblock {\em Stat. Med.\/}~{\em 28}, 1108--1130.

\bibitem[\protect\citeauthoryear{Gilbert and Hudgens}{Gilbert and
  Hudgens}{2008}]{gilbert2008evaluating}
Gilbert, P.~B. and M.~G. Hudgens (2008).
\newblock Evaluating candidate principal surrogate endpoints.
\newblock {\em Biometrics\/}~{\em 64}, 1146--1154.

\bibitem[\protect\citeauthoryear{Goodman}{Goodman}{1974}]{goodman1974explorato%
ry}
Goodman, L.~A. (1974).
\newblock Exploratory latent structure analysis using both identifiable and
  unidentifiable models.
\newblock {\em Biometrika\/}~{\em 61}, 215--231.

\bibitem[\protect\citeauthoryear{Gustafson}{Gustafson}{2009}]{gustafson2009lim%
its}
Gustafson, P. (2009).
\newblock What are the limits of posterior distributions arising from
  nonidentified models, and why should we care?
\newblock {\em J. Am. Statist. Ass.\/}~{\em 104}, 1682--1695.

\bibitem[\protect\citeauthoryear{Huang and Gilbert}{Huang and
  Gilbert}{2011}]{huang2011comparing}
Huang, Y. and P.~B. Gilbert (2011).
\newblock Comparing biomarkers as principal surrogate endpoints.
\newblock {\em Biometrics\/}~{\em 67}, 1442--1451.

\bibitem[\protect\citeauthoryear{Joffe and Greene}{Joffe and
  Greene}{2009}]{joffe2009related}
Joffe, M.~M. and T.~Greene (2009).
\newblock Related causal frameworks for surrogate outcomes.
\newblock {\em Biometrics\/}~{\em 65}, 530--538.

\bibitem[\protect\citeauthoryear{Ju and Geng}{Ju and
  Geng}{2010}]{ju2010criteria}
Ju, C. and Z.~Geng (2010).
\newblock Criteria for surrogate end points based on causal distributions.
\newblock {\em J. R. Statist. Soc. B\/}~{\em 72}, 129--142.

\bibitem[\protect\citeauthoryear{Lauritzen}{Lauritzen}{2004}]{lauritzen2004dis%
cussion}
Lauritzen, S.~L. (2004).
\newblock Discussion on causality.
\newblock {\em Scand. J. Stat.\/}~{\em 31}, 189--193.

\bibitem[\protect\citeauthoryear{Li, Taylor, Elliott, and Sargent}{Li
  et~al.}{2011}]{li2011causal}
Li, Y., J.~M. Taylor, M.~R. Elliott, and D.~J. Sargent (2011).
\newblock Causal assessment of surrogacy in a meta-analysis of colorectal
  cancer trials.
\newblock {\em Biostatistics\/}~{\em 12}, 478--492.

\bibitem[\protect\citeauthoryear{Liu}{Liu}{2001}]{liu2001generalized}
Liu, J.~S. (2001).
\newblock {\em Monte Carlo Strategies in Scientific Computing}.
\newblock Springer: New York.

\bibitem[\protect\citeauthoryear{Mattei, Li, and Mealli}{Mattei
  et~al.}{2013}]{mattei2012exploiting}
Mattei, A., F.~Li, and F.~Mealli (2013).
\newblock Exploiting multiple outcomes in {B}ayesian inference for causal
  effects with intermediate variables.
\newblock {\em Ann. Appl. Stat.\/}~{\em 7}, 2336--2360.

\bibitem[\protect\citeauthoryear{Mattei and Mealli}{Mattei and
  Mealli}{2011}]{mattei2011augmented}
Mattei, A. and F.~Mealli (2011).
\newblock Augmented designs to assess principal strata direct effects.
\newblock {\em J. R. Statist. Soc. B\/}~{\em 73}, 729--752.

\bibitem[\protect\citeauthoryear{Mattei, Mealli, and Pacini}{Mattei
  et~al.}{2014}]{mattei2014identification}
Mattei, A., F.~Mealli, and B.~Pacini (2014).
\newblock Identification of causal effects in the presence of nonignorable
  missing outcome values.
\newblock {\em Biometrics\/}~{\em 70}, 278--288.

\bibitem[\protect\citeauthoryear{Mealli and Mattei}{Mealli and
  Mattei}{2012}]{mealli2012refreshing}
Mealli, F. and A.~Mattei (2012).
\newblock A refreshing account of principal stratification.
\newblock {\em Int. J. Biostatistics\/}~{\em 8}, 1--37.

\bibitem[\protect\citeauthoryear{Mealli and Pacini}{Mealli and
  Pacini}{2008}]{mealli2008comparing}
Mealli, F. and B.~Pacini (2008).
\newblock Comparing principal stratification and selection models in parametric
  causal inference with nonignorable missingness.
\newblock {\em Compute. Stat. Data An.\/}~{\em 53}, 507--516.

\bibitem[\protect\citeauthoryear{Mealli and Pacini}{Mealli and
  Pacini}{2013}]{mealli2013using}
Mealli, F. and B.~Pacini (2013).
\newblock Using secondary outcomes to sharpen inference in randomized
  experiments with noncompliance.
\newblock {\em J. Am. Statist. Ass.\/}~{\em 108}, 1120--1131.

\bibitem[\protect\citeauthoryear{Meng}{Meng}{1994}]{meng1994posterior}
Meng, X.~L. (1994).
\newblock Posterior predictive $p$-values.
\newblock {\em Ann. Stat.\/}~{\em 22}, 1142--1160.

\bibitem[\protect\citeauthoryear{Prentice}{Prentice}{1989}]{prentice1989surrog%
ate}
Prentice, R.~L. (1989).
\newblock Surrogate endpoints in clinical trials: definition and operational
  criteria.
\newblock {\em Stat. Med.\/}~{\em 8}, 431--440.

\bibitem[\protect\citeauthoryear{Qin, Gilbert, Follmann, and Li}{Qin
  et~al.}{2008}]{qin2008assessing}
Qin, L., P.~B. Gilbert, D.~Follmann, and D.~Li (2008).
\newblock Assessing surrogate endpoints in vaccine trials with case-cohort
  sampling and the {C}ox model.
\newblock {\em Ann. Appl. Stat\/}~{\em 2}, 386.

\bibitem[\protect\citeauthoryear{Roche, Miglioretti, Zeger, and Rathouz}{Roche
  et~al.}{1997}]{bandeen1997latent}
Roche, K., D.~L. Miglioretti, S.~L. Zeger, and P.~J. Rathouz (1997).
\newblock Latent variable regression for multiple discrete outcomes.
\newblock {\em J. Am. Statist. Ass.\/}~{\em 92}, 1375--1386.

\bibitem[\protect\citeauthoryear{Romano and Shaikh}{Romano and
  Shaikh}{2010}]{romano2010inference}
Romano, J.~P. and A.~M. Shaikh (2010).
\newblock Inference for the identified set in partially identified econometric
  models.
\newblock {\em Econometrica\/}~{\em 78}, 169--211.

\bibitem[\protect\citeauthoryear{Rubin}{Rubin}{1980}]{rubin1980randomization}
Rubin, D.~B. (1980).
\newblock Comment on ``{R}andomization analysis of experimental data: The
  {F}isher randomization test''.
\newblock {\em J. Am. Statist. Ass.\/}~{\em 75}, 591--593.

\bibitem[\protect\citeauthoryear{Rubin}{Rubin}{1984}]{rubin1984bayesianly}
Rubin, D.~B. (1984).
\newblock Bayesianly justifiable and relevant frequency calculations for the
  applied statistician.
\newblock {\em Ann. Stat.\/}~{\em 12}, 1151--1172.

\bibitem[\protect\citeauthoryear{Rubin}{Rubin}{2004}]{rubin2004direct}
Rubin, D.~B. (2004).
\newblock Direct and indirect causal effects via potential outcomes (with
  discussion).
\newblock {\em Scand. J. Stat.\/}~{\em 31}, 161--170.

\bibitem[\protect\citeauthoryear{Rubin}{Rubin}{2006}]{rubin2006causal}
Rubin, D.~B. (2006).
\newblock Causal inference through potential outcomes and principal
  stratification: application to studies with ``censoring" due to death (with
  discussion).
\newblock {\em Stat. Sci.\/}~{\em 21}, 299--309.

\bibitem[\protect\citeauthoryear{Sargent, Wieand, Haller, et~al.}{Sargent
  et~al.}{2005}]{sargent2005disease}
Sargent, D.~J., H.~S. Wieand, D.~G. Haller, et~al. (2005).
\newblock Disease-free survival versus overall survival as a primary end point
  for adjuvant colon cancer studies: individual patient data from 20,898
  patients on 18 randomized trials.
\newblock {\em J. Clin. Oncol.\/}~{\em 23}, 8664--8670.

\bibitem[\protect\citeauthoryear{Schwartz, Li, and Mealli}{Schwartz
  et~al.}{2011}]{schwartz2011bayesian}
Schwartz, S.~L., F.~Li, and F.~Mealli (2011).
\newblock A {B}ayesian semiparametric approach to intermediate variables in
  causal inference.
\newblock {\em J. Am. Statist. Ass.\/}~{\em 106}, 1331--1344.

\bibitem[\protect\citeauthoryear{Skrondal and Rabe-Hesketh}{Skrondal and
  Rabe-Hesketh}{2004}]{skrondal2004generalized}
Skrondal, A. and S.~Rabe-Hesketh (2004).
\newblock {\em Generalized Latent Variable Modeling: Multilevel, Longitudinal,
  and Structural Equation Models}.
\newblock CRC Press: London.

\bibitem[\protect\citeauthoryear{VanderWeele}{VanderWeele}{2013}]{vanderweele2%
013surrogate}
VanderWeele, T.~J. (2013).
\newblock Surrogate measures and consistent surrogates (with discussion).
\newblock {\em Biometrics\/}~{\em 69}, 561--565.

\bibitem[\protect\citeauthoryear{Yang and Small}{Yang and
  Small}{2015}]{yang2015using}
Yang, F. and D.~S. Small (2015).
\newblock Using post-quality of life measurement information in censoring by
  death problems.
\newblock {\em J. R. Statist. Soc. B\/}, in press.

\bibitem[\protect\citeauthoryear{Zhang and Rubin}{Zhang and
  Rubin}{2003}]{zhang2003estimation}
Zhang, J.~L. and D.~B. Rubin (2003).
\newblock Estimation of causal effects via principal stratification when some
  outcomes are truncated by ``death''.
\newblock {\em J. Educ. Behav. Stat.\/}~{\em 28}, 353--368.

\bibitem[\protect\citeauthoryear{Zhang, Rubin, and Mealli}{Zhang
  et~al.}{2009}]{zhang2009likelihood}
Zhang, J.~L., D.~B. Rubin, and F.~Mealli (2009).
\newblock Likelihood-based analysis of causal effects of job-training programs
  using principal stratification.
\newblock {\em J. Am. Statist. Ass.\/}~{\em 104}, 166--176.

\bibitem[\protect\citeauthoryear{Zigler and Belin}{Zigler and
  Belin}{2012}]{zigler2012bayesian}
Zigler, C.~M. and T.~R. Belin (2012).
\newblock A {B}ayesian approach to improved estimation of causal effect
  predictiveness for a principal surrogate endpoint.
\newblock {\em Biometrics\/}~{\em 68}, 922--932.

\end{thebibliography}

\newpage 
{\bf \Huge Supplementary Materials}

\bigskip

Appendix A  derives the bounds for the PSACEs both with and without monotonicity, and shows the application of these bounds to the ACCTs data. Appendix B provides the proofs for the propositions and theorems. Appendix C gives the computational details. Appendix D provides more details about the simulation studies.

 We re-introduce all the notation used in the main text as follows:
\begin{alignat*}{3}
 P_{zsr}&= P(S=s \mid Z=z, R=r),  && \quad Q_{zsr}=P(Y=1 \mid Z=z, S=s, R=r),\\
  \omega_{ys \mid z r}&= P(Y=y, S=s\mid Z=z, R=r), &&\quad p_r = P(R=r), \\
  \alpha_r &= P(Z=1\mid R=r), &&\quad \pi_{ur}=P(U=u\mid R=r),\\
    \delta_{zur} &= P(Y=1\mid Z=z, U=u, R=r), &&\quad ACE_{ur} = E\{ Y(1)-Y(0)\mid  U=u, R=r\}. 
\end{alignat*}
Under the homogeneity assumption, we have: 
$$\delta_{zu} = P(Y=1\mid Z=z, U=u)= \delta_{zur}, \quad ACE_{u} = E\{ Y(1)-Y(0)\mid  U=u\}=ACE_{ur}.
$$
We further define $\bm{p} = \{p_r: r=1,\ldots,N_R\}$, $\bm{\alpha}=\{\alpha_r: r=1,\ldots,N_R\}$,
$\bm{\pi}_r=\{\pi_{ss,r}, \pi_{s\bar{s},r}, \pi_{\bar{s}\bar{s},r}, \pi_{\bar{s}s,r} \}$,
$\bm{\pi}=\{  \bm{\pi}_{r}: r=1,\ldots,N_R\}$,
and $ \bm{\delta}=\{\delta_{zu }: z=0,1; u=ss,s\bar{s},\bar{s}\bar{s}, \bar{s}s  \}$.

\section*{\bf Appendix A: Bounds for PSACEs}
We need the following lemma to simplify the derivations of the bounds.
\begin{lemma}
\label{lem:1}
Let $X_0$ be a mixture of two Bernoulli distributions $X_1$ and $X_2$ with $X_0\sim \alpha X_1+(1-\alpha)X_2, X_i\sim $ Bernoulli$(p_i)$, and a known mixing proportion $\alpha$. 
Then we have
\begin{eqnarray*}
&&\max\left( 0,1-\frac{1-p_0}{\alpha}\right) \leq p_1\leq \min\left(1,\frac{p_0}{\alpha}\right),\\
&&\max\left( 0,1-\frac{1-p_0}{1-\alpha}\right) \leq p_2\leq \min\left(1,\frac{p_0}{1-\alpha}\right).
\end{eqnarray*}
\end{lemma}

\noindent {\it Proof of Lemma 1.} See \citet{cheng2006bounds}. \qed

\renewcommand {\theequation} {A.\arabic{equation}}

\subsection*{\bf Appendix A.1: Bounds without monotonicity}
\renewcommand {\theproposition} {A.\arabic{proposition}}
Let $\text{Upper}_{nm}(ACE_{ur})$ and $\text{Lower}_{nm}(ACE_{ur})$  denote the upper and lower bounds of $ACE_{ur}$ without monotonicity, respectively.
\begin{proposition} \label{prop::bounds1}
Without monotonicity, the bounds of $ACE_{ss,r}$ are
\begin{eqnarray*}
\text{Upper}_{nm}(ACE_{ss,r}) &=&
\min\left\{1,\frac{Q_{11r}   P_{11r}}{P_{01r}-P_{10r}}\right\}+\min\left\{0,\frac{(1-Q_{01r})   P_{01r}}{P_{01r}-P_{10r}}-1\right\}, \\
\text{Lower}_{nm}(ACE_{ss,r}) &=&
\max\left\{0,1-\frac{(1-Q_{11r})  P_{11r}}{P_{01r}-P_{10r}}\right\}-\min\left\{1,\frac{Q_{01r}  
P_{01}}{P_{01r}-P_{10r}}\right\},
\end{eqnarray*}
if $P_{01r} \geq P_{10r}$, and $[-1,1]$ if $P_{01r} < P_{10r}$.
The bounds of $ACE_{s\bar{s},r}$ are
\begin{eqnarray*}
\text{Upper}_{nm}(ACE_{s\bar{s},r}) &=&
\min\left\{1,\frac{Q_{11r}   P_{11r}}{P_{11r}-P_{01r}}\right\}-\max\left\{0,\frac{Q_{00r}   P_{00r}-P_{10r}}{P_{11r}-P_{01r}}\right\},\\
\text{Lower}_{nm}(ACE_{s\bar{s},r}) &=&
\max\left\{0,\frac{Q_{11r}P_{11r}}{P_{11r}-P_{01r}}-1\right\}-\min\left\{1,\frac{Q_{00r}   P_{00r}}{P_{11r}-P_{01r}}\right\},
\end{eqnarray*}
if $P_{01r}\leq P_{11r}$, and $[-1,1]$ if $P_{01r} > P_{11r}$.
The bounds of $ACE_{\bar{s}\bar{s},r}$ are:
\begin{eqnarray*}
\text{Upper}_{nm}(ACE_{\bar{s}\bar{s},r}) &=&
\min\left\{1,\frac{Q_{10r}   P_{10r}}{P_{10r}-P_{01r}}\right\}-\max\left\{0,\frac{Q_{00r}   P_{00r}-P_{11r}}{P_{10r}-P_{01r}}\right\}, \\
\text{Lower}_{nm}(ACE_{\bar{s}\bar{s},r}) &=&
\max\left\{0,\frac{Q_{10r}P_{10r}-P_{01r}}{P_{10r}-P_{01r}}\right\}-\min\left\{1,\frac{Q_{00r}   P_{00r}}{P_{10r}-P_{01r}}\right\},
\end{eqnarray*}
if $P_{01r}\leq P_{10r}$, and $[-1,1]$ if $P_{01r} > P_{10r}$.
The bounds of $ACE_{\bar{s}s,r}$ are:
\begin{eqnarray*}
\text{Upper}_{nm}(ACE_{\bar{s}s,r}) &=&
\min\left\{1,\frac{Q_{10r}   P_{10r}}{P_{01r}-P_{11r}}\right\}-\max\left\{0,\frac{Q_{01r}   P_{01r}-P_{11r}}{P_{01r}-P_{11r}}\right\},\\
\text{Lower}_{nm}(ACE_{\bar{s}s,r}) &=&
\max\left\{0,\frac{Q_{10r}P_{10r}-P_{00r}}{P_{01r}-P_{11r}}\right\}-\min\left\{1,\frac{Q_{01r}   P_{01r}}{P_{01r}-P_{11r}}\right\},
\end{eqnarray*}
if $P_{01r}\geq P_{11r}$, and $[-1,1]$ if $P_{01r} < P_{11r}$.
\end{proposition}

\noindent {\it Proof of Proposition \ref{prop::bounds1}.}
Without monotonicity, each of the subpopulation with $(Z=z,S=s)$ is a mixture of two latent principal strata, implying that
\begin{eqnarray}
 \pi_{\bar{s}\bar{s},r}+\pi_{s\bar{s},r}&=&P_{00r},  \label{eq::prop1-append}\\
\pi_{ss,r}+\pi_{\bar{s}s,r}&=&P_{01r},\\
 \pi_{\bar{s}\bar{s},r}+\pi_{\bar{s}s,r}&=&P_{10r}, \\
  \pi_{ss,r}+\pi_{s\bar{s},r}&=&P_{11r},\\
 \pi_{ss,r}+\pi_{s\bar{s},r}+\pi_{\bar{s}\bar{s},r}+\pi_{\bar{s}s,r}&=&1, \label{eq::prop2-append}\\
\delta_{1,ss,r}   \frac{\pi_{ss,r}}{P_{11r}}+\delta_{1,s\bar{s},r}  \frac{\pi_{s\bar{s},r}}{P_{11r}} &=& Q_{11r}, \label{eqn:A1} \\
\delta_{0,ss,r}   \frac{\pi_{ss,r}}{P_{01r}}+\delta_{0,\bar{s}s,r}  \frac{\pi_{\bar{s}s,r}}{P_{01r}} &=& Q_{01r}. \label{eqn:A2}
\end{eqnarray}
Supposing that we know $\pi_{\bar{s}s,r}$ in Equations (\ref{eq::prop1-append}) and (\ref{eq::prop2-append}), we can express $\pi_{ur}$'s in terms of $\pi_{\bar{s}s,r}$ :
\begin{eqnarray}
 \pi_{ss,r}&=&P_{01r}-\pi_{\bar{s}s,r},  \label{eq::prob1}  \\
 \pi_{\bar{s}\bar{s},r}&=&P_{10r}-\pi_{\bar{s}s,r}, \\
 \pi_{s\bar{s},r}&=&1-(\pi_{ss,r}+\pi_{\bar{s}\bar{s},r}+\pi_{\bar{s}s,r})=P_{11r}-P_{01r}+\pi_{\bar{s}s,r}. \label{eq::prob2}
\end{eqnarray}
Since $0 \leq \pi_{ur}\leq 1$, we
can obtain the bounds for $\pi_{\bar{s}s,r}$ from Equations (\ref{eq::prob1}) to (\ref{eq::prob2}):
$$ 
\pi_{\bar{s}s, r}  \in \left[   \max ( 0,P_{01r}-P_{11r} ) ,  \min (P_{01r},P_{10r} )  \right]  \equiv  \mathcal{R}_{r}.
$$
From Lemma \ref{lem:1} and Equations (\ref{eqn:A1}) and (\ref{eqn:A2}), we can get the bounds of  $ACE_{ss,r}$ for a given $\pi_{\bar{s}s,r}$:
\begin{eqnarray*}
&&\text{Upper}_{nm}(ACE_{ss,r}\mid  \pi_{\bar{s}s, r}) \\
&=&\min\left\{1,\frac{Q_{11r}   P_{11r}}{P_{01r}-\pi_{\bar{s}s,r}}\right\}-\max\left\{0,\frac{Q_{01r}   P_{01r}-\pi_{\bar{s}s,r}}{P_{01r}-\pi_{\bar{s}s,r}}\right\}\\
&=&\min\left\{1,\frac{Q_{11r}   P_{11r}}{P_{01r}-\pi_{\bar{s}s,r}}\right\}+\min\left\{0,\frac{(1-Q_{01r})  P_{01r}}{P_{01r}-\pi_{\bar{s}s,r}}-1\right\},
\end{eqnarray*}
\begin{eqnarray*}
&&\text{Lower}_{nm}(ACE_{ss,r}\mid \pi_{\bar{s}s, r})\\
&=&\max\left\{0,\frac{Q_{11r}P_{11r}-P_{11r}+P_{01r}-\pi_{\bar{s}s,r}}{P_{01r}-\pi_{\bar{s}s,r}}\right\}-\min\left\{1,\frac{Q_{01r}  
P_{01r}}{P_{01r}-\pi_{\bar{s}s,r}}\right\}\\
&=&\max\left\{0,1-\frac{(1-Q_{11r})  P_{11r}}{P_{01r}-\pi_{\bar{s}s,r}}\right\}-\min\left\{1,\frac{Q_{01r}  
P_{01r}}{P_{01r}-\pi_{\bar{s}s,r}}\right\}.
\end{eqnarray*}
The bounds for $ACE_{ss,r}$ can be obtained by maximizing or minimizing the above bounds over $\mathcal{R}_r$, the feasible region of $\pi_{\bar{s}s, r}$, namely,
\begin{eqnarray*}
&&\text{Upper}_{nm}(ACE_{ss,r}) \\
&=&\max_{ \pi_{\bar{s}s, r}  \in \mathcal{R}_{r} }
\Bigg[\min\left\{1,\frac{Q_{11r}   P_{11r}}{P_{01r}-\pi_{\bar{s}s,r}}\right\}+\min\left\{0,\frac{(1-Q_{01r})  P_{01r}}{P_{01r}-\pi_{\bar{s}s,r}}-1\right\}\Bigg],\\
&&\text{Lower}_{nm}(ACE_{ss,r}) \\
&=&\min_{ \pi_{\bar{s}s, r}  \in  \mathcal{R}_{r}  }
\Bigg[\max\left\{0,1-\frac{(1-Q_{11r}) P_{11r}}{P_{01r}-\pi_{\bar{s}s,r}}\right\}-\min\left\{1,\frac{Q_{01r}  
P_{01r}}{P_{01r}-\pi_{\bar{s}s,r}}\right\}\Bigg].
\end{eqnarray*}
Since 
$$
\min \Big \{1,\frac{Q_{11r}  P_{11r}}{P_{01r}-\pi_{\bar{s}s,r}}\Big\}+\min\Big \{0,\frac{(1-Q_{01r})  P_{01r}}{P_{01r}-\pi_{\bar{s}s,r}}-1\Big \}
$$ 
is increasing in $\pi_{\bar{s}s,r}$, and 
$$
\max\Big\{0,1-\frac{(1-Q_{11r})  P_{11r}}{P_{01r}-\pi_{\bar{s}s,r}}\Big\}-\min\Big\{1,\frac{Q_{01r}  
P_{01r}}{P_{01r}-\pi_{\bar{s}s,r}}\Big\}
$$ 
is decreasing in $\pi_{\bar{s}s,r}$, the above bounds can be simplified as
\begin{eqnarray*}
\text{Upper}_{nm}(ACE_{ss,r}) &=&
\min\left\{1,\frac{Q_{11r}   P_{11r}}{P_{01r}-P_{10r}}\right\}+\min\left\{0,\frac{(1-Q_{01r})   P_{01r}}{P_{01r}-P_{10r}}-1\right\}, \\
\text{Lower}_{nm}(ACE_{ss,r}) &=&
\max\left\{0,1-\frac{(1-Q_{11r})  P_{11r}}{P_{01r}-P_{10r}}\right\}-\min\left\{1,\frac{Q_{01r}  
P_{01}}{P_{01r}-P_{10r}}\right\},
\end{eqnarray*}
if $P_{01r} \geq P_{10r}$, and $[-1,1]$ if $P_{01r} < P_{10r}$.

Similarly, we can obtain the bounds of the other three principal strata.\qed

\subsection*{\bf Appendix A.2: Bounds with monotonicity}
\renewcommand {\theproposition} {A.2}

Let $\text{Upper}_{m}(ACE_{ur})$ and $\text{Lower}_{m}(ACE_{ur})$  denote the upper and lower bounds of $ACE_{ur}$ without monotonicity, respectively.
\begin{proposition} \label{prop::bounds2}
With monotonicity, the bounds of $ACE_{ss,r}$ are
\begin{eqnarray*}
&&\text{Upper}_m(ACE_{ss,r})=\min\left\{1,\frac{Q_{11r}   P_{11r}}{P_{01r}}\right\}-Q_{01r},\\
&&\text{Lower}_m(ACE_{ss,r})=\max\left\{0,\frac{Q_{11r}P_{11r}-P_{11r}+P_{01r}}{P_{01r}}\right\}-Q_{01r};
\end{eqnarray*}
the bounds for $ACE_{s\bar{s},r}$ are
\begin{eqnarray*}
&&\text{Upper}_{m}(ACE_{s\bar{s},r})=\min\left\{1,\frac{Q_{11r}   P_{11r}}{P_{11r}-P_{01r}}\right\}-\max\left\{0,\frac{Q_{00r}   P_{00r}-P_{10r}}{P_{11r}-P_{01r}}\right\},\\
&&\text{Lower}_{m}(ACE_{s\bar{s},r})=\max\left\{0,\frac{Q_{11r}P_{11}-P_{10r}}{P_{01r}}\right\}-\min\left\{1,\frac{Q_{00r}   P_{00r}}{P_{11}-P_{01r}}\right\};
\end{eqnarray*}
and the bounds for $ACE_{ss,r}$ are
\begin{eqnarray*}
&&\text{Upper}_{m}(ACE_{\bar{s}\bar{s},r})=Q_{10r} -\max\left\{0,\frac{Q_{00r}   P_{00r}-P_{11r}+P_{10r}}{P_{10r}}\right\},\\
&&\text{Lower}_{m}(ACE_{\bar{s}\bar{s},r})=Q_{10r}-\min\left\{1,\frac{Q_{00r}   P_{00r}}{P_{10r}}\right\}.
\end{eqnarray*}
\end{proposition}

\noindent {\it Proof of Proposition \ref{prop::bounds2}.}
When monotonicity holds, we have $\pi_{\bar{s}s,r}=0$, and we can identify the proportions of all the principal strata,  and thus we have 
\begin{eqnarray*}
\text{Upper}_m(ACE_{ss,r})=\text{Upper}_{nm}(ACE_{ss,r}\mid  \pi_{\bar{s}s, r}=0),\\
\text{Lower}_m(ACE_{ss,r})=\text{Lower}_{nm}(ACE_{ss,r}\mid  \pi_{\bar{s}s, r}=0).
\end{eqnarray*}
Similarly, we can obtain the bounds of the other two principal strata.\qed

\subsection*{\bf Appendix A.3: Bounds for the ACCTs data}

We compute bounds of the PSACEs for the $10$ trials separately. 
Figure \ref{fig:2} shows the large sample bounds and confidence intervals for PSACEs based on the bootstrap \citep{cheng2006bounds}.
Figure \ref{fig:2a} shows that, without monotonicity, bounds of the PSACEs are barely informative. These bounds always contain zero, and bounds of $ACE_{\bar{s}\bar{s},r}$ and $ACE_{\bar{s}s,r}$ are $[-1,1]$, which are too wide for us to get any useful information.
Figure \ref{fig:2b} shows that, with monotonicity, bounds are sometimes informative. But the confidence intervals for the PSACEs always contain zero, which do not provide strong evidence for the presence of causal effects. 

Standard resampling techniques, including the bootstrap, may be inconsistent and lead to confidence sets that are not asymptotically valid in a point-wise or uniform sense for drawing inference on partially identified quantities \citep{andrews2000inconsistency,andrews2009validity,romano2010inference}. 
Here we follow \citet{cheng2006bounds}'s procedure to obtain the confidence intervals for the bounds, which is easy to implement.  For more careful analysis, we can use the methods developed by \citet{chernozhukov2013intersection}. However, we do not do this here because calculating the bounds is not our primary goal. And what is more, the bounds are barely informative in our application.

\section*{\bf Appendix B: Proofs of the propositions and theorems}
\setcounter{equation}{0}
\renewcommand {\theequation} {B.\arabic{equation}}
\setcounter{proposition}{0}
\renewcommand {\theproposition} {\arabic{proposition}}

\subsection*{\bf Appendix B.1: Proof of Proposition 1 in Section 3}
\begin{proposition}
\label{prop:2}
Under Assumption 1, Assumption 3 is equivalent to $R \ind Y \mid(U, Z)$, which implies $ACE_{ur} = ACE_{ur'}$ for all $u=ss$,
$s\bar{s}$, $\bar{s}s$, $\bar{s}\bar{s}$ and $r \neq r'$.
\end{proposition}

\noindent {\it Proof of Proposition \ref{prop:2}.} Under randomization and homogeneity (Assumption 3), we have 
\begin{eqnarray*}
&&P(Y=1\mid U=u,R=r,Z=z) \\
&=& P\{Y(z)=1\mid U=u,R=r,Z=z\}\\
&=& P\{Y(z)=1\mid U=u,R=r\} = P\{Y(z)=1\mid U=u\},
\end{eqnarray*}
and
\begin{eqnarray*}
&&P(Y=1\mid U=u,Z=z )\\
&=&\sum _{r}P(Y=1\mid U=u,R=r,Z=z)\cdot P(R=r\mid U=u,Z=z)\\
&=&\sum _{r}P\{Y(z)=1\mid U=u\}\cdot P(R=r\mid U=u,Z=z)\\
&=&P\{Y(z)=1\mid U=u\}\cdot \sum _{r} P(R=r\mid U=u,Z=z)\\
&=&P\{Y(z)=1\mid U=u\}=P(Y=1\mid U=u,R=r,Z=z),
\end{eqnarray*}
implying $R \ind Y \mid (U,Z ) $.
Thus we have 
$$P(Y=1\mid U=u,R=r,Z=z)=P(Y=1\mid U=u,R=r',Z=z)$$ for $z=0,1$; $u=ss,
s\bar{s}, \bar{s}s, \bar{s}\bar{s}$; and $r,r'=1,\ldots,N_R$.
Under homogeneity (Assumption 3), we have $E\{Y(z)\mid U=u,R=r\} = E\{Y(z)\mid U=u\} $, and thus $ACE_{ur}=ACE_{ur'}$, for $z=0,1$; $u=ss,
s\bar{s}, \bar{s}s, \bar{s}\bar{s}$; and $r,r'=1,\ldots,N_R$.
In addition, under randomization and homogeneity (Assumption 3), we have 
$$
P\{Y(z)\mid U=u,R=r\}=P(Y\mid Z=z,U=u,R=r)=P(Y\mid Z=z,U=u),
$$
and thus $ACE_{u}=\delta_{1u}-\delta_{0u}$ does not depend on $r.$\qed

\subsection*{\bf Appendix B.2: Proof of Theorem 1 in Section 3.1}
\begin{thm}
\label{thm:1}
Under Assumptions 1 to 4, we have that, for $u=ss$, $s\bar{s}$ and $\bar{s}\bar{s}$,
\begin{enumerate}
[(a)]
\item
$P(Y=1\mid Z=1,U=u)$ for all principal strata $u$ are identifiable if Assumption 4(a) holds;

\item
$P(Y=1\mid Z=0,U=u)$ for all principal strata $u$ are identifiable if Assumption 4(b) holds; 

\item
$ACE_{u}$ for all principal strata $u$ are identifiable if both Assumptions 4(a) and 4(b) hold.
\end{enumerate}
\end{thm}

\noindent {\it Proof of Theorem \ref{thm:1}.} First, we can identify the proportions of all the principal strata under monotonicity.
Let  $\omega_{ys\mid zr}= P(Y=y,S=s\mid Z=z,R=r)$ denote the probability of $(Y=y, S=s)$ given $(Z=z, R=r)$, which can be identified from the observed data.

The subpopulation with $(Z=1, S=0)$ is equivalent to stratum $\bar{s}\bar{s}$ under treatment, and the subpopulation with $(Z=0, S=1)$ is equivalent to stratum $ss$ under control. Homogeneity (Assumption 3) implies that
$
\omega_{10\mid 1r}=\delta_{1,\bar{s}\bar{s}}\cdot \pi_{\bar{s}\bar{s},r}$
and
$
\omega_{11\mid 0r}=\delta_{0,ss}\cdot \pi_{ss,r},
$
and therefore we can identify $\delta_{1,\bar{s}\bar{s}} = \omega_{10\mid 1r}/\pi_{\bar{s}\bar{s},r}$ and $ \delta_{0,ss} = \omega_{11\mid 0r}/\pi_{ss,r}$ from the observed data directly.

The subpopulation with $(Z=1, S=1)$ is a mixture of two latent strata $ss$ and $s\bar{s}$, and under homogeneity (Assumption 3) we have:
\begin{eqnarray*}
\omega_{11\mid 1r}
= \delta_{1,ss}\cdot \pi_{ss,r}+\delta_{1,s\bar{s}}\cdot \pi_{s\bar{s},r}.
\end{eqnarray*}
For the two trials $r_1$ and $r_2$ satisfying Assumption 4(a) with $\pi_{ss,r_1}\pi_{s\bar{s},r_2}\neq \pi_{s\bar{s},r_1}\pi_{ss,r_2}$, we have:
$$\left\{ 
\begin{array}{ccc}
\omega_{11\mid 1r_1}&=&\delta_{1,ss}\cdot \pi_{ss,r_1}+\delta_{1,s\bar{s}}\cdot \pi_{s\bar{s},r_1},\\
\omega_{11\mid 1r_2}&=&\delta_{1,ss}\cdot \pi_{ss,r_2}+\delta_{1,s\bar{s}}\cdot \pi_{s\bar{s},r_2},
\end{array} 
\right.
$$
from which we can find a unique solution for $(\delta_{1,ss}, \delta_{1,s\bar{s}})$:
$$
\left\{
\begin{array}{ccc}
\delta_{1,ss} &=&   \frac{\pi_{s\bar{s},r_2}\omega_{11\mid 1r_1}-\pi_{s\bar{s},r_1}\omega_{11\mid 1r_2}}{\pi_{ss,r_1}\pi_{s\bar{s},r_2}- \pi_{s\bar{s},r_1}\pi_{ss,r_2}},\\
\delta_{1,s\bar{s}} &=&   \frac{\pi_{ss,r_1}\omega_{11\mid 1r_2}-\pi_{ss,r_2}\omega_{11\mid 1r_1}}{\pi_{ss,r_1}\pi_{s\bar{s},r_2}- \pi_{s\bar{s},r_1}\pi_{ss,r_2}}.
\end{array}
\right.
$$
Similarly, the subpopulation with $(Z=0, S=0)$ is a mixture of two latent strata $s\bar{s}$ and $\bar{s}\bar{s}$, and under homogeneity (Assumption 3) we have:
\begin{eqnarray*}
\omega_{10\mid 0r}
= \delta_{0,\bar{s}\bar{s}}\cdot \pi_{\bar{s}\bar{s},r}+\delta_{0,s\bar{s}}\cdot \pi_{s\bar{s},r}.
\end{eqnarray*}
For the two trials $r_3$ and $r_4$ satisfying Assumption 4(b) with $\pi_{\bar{s}\bar{s},r_3}\pi_{s\bar{s},r_4} \neq \pi_{s\bar{s},r_3}\pi_{\bar{s}\bar{s},r_4}$, we have:
$$
\left\{
\begin{array}{ccc}
\omega_{10\mid 0r_3} &=& \delta_{0,\bar{s}\bar{s}}\cdot \pi_{\bar{s}\bar{s},r_3}+\delta_{0,s\bar{s}}\cdot \pi_{s\bar{s},r_3},\\
\omega_{10\mid 0r_4} &=& \delta_{0,\bar{s}\bar{s}}\cdot \pi_{\bar{s}\bar{s},r_4}+\delta_{0,s\bar{s}}\cdot \pi_{s\bar{s},r_4},
\end{array}
\right.
$$
from which we can obtain a unique solution to $(\delta_{0,\bar{s}\bar{s}}, \delta_{0,s\bar{s}})$:
$$
\left\{
\begin{array}{ccc}
\delta_{0,\bar{s}\bar{s}} &=&   \frac{\pi_{s\bar{s},r_4}\omega_{10\mid 0r_3}-\pi_{s\bar{s},r_3}\omega_{10\mid 0r_4}}{\pi_{\bar{s}\bar{s},r_3}\pi_{s\bar{s},r_4}- \pi_{s\bar{s},r_3}\pi_{\bar{s}\bar{s},r_4}},\\
\delta_{0,s\bar{s}} &=&   \frac{\pi_{\bar{s}\bar{s},r_3}\omega_{10\mid 0r_4}-\pi_{\bar{s}\bar{s},r_4}\omega_{10\mid 0r_3}}{\pi_{\bar{s}\bar{s},r_3}\pi_{s\bar{s},r_4}- \pi_{s\bar{s},r_3}\pi_{\bar{s}\bar{s},r_4}}.
\end{array}
\right.
$$
Since $ACE_{u}=\delta_{1u}-\delta_{0u}$, we can identify $ACE_{u}$'s under Assumptions 1 to 4.  \qed

\subsection*{\bf Appendix B.3: Proof of Proposition 2 in Section 3.2}
\renewcommand {\thetable} {B.3.\arabic{table}}
\renewcommand {\theexample} {B.3.\arabic{example}}
\begin{proposition}
\label{prop:nec}
Under Assumptions 1 and 3,
a necessary condition for the local identifiability of the joint distribution of  $(Y,Z,S,U,R)$ is $N_R\geq 3$.
\end{proposition}

\noindent {\it Proof of Proposition \ref{prop:nec}.} The observed data of $(R, Z, S, Y)$ form an $N_R\times 2\times 2\times 2$ contingency table with $(8N_R-1)$ free frequencies, therefore the saturated model contains $(8N_R - 1)$ parameters. The joint distributions contains $(5N_R + 7)$ parameters, with $P(R)$ contributing $(N_R-1)$ parameters, $P(Z\mid R)$ contributing $N_R$ parameters, $P(U\mid R)$ contributing $3N_R$ parameters, and $P(Y\mid Z,U)$ contributing $8$ parameters.
Therefore, the model with monotonicity has $  (8N_R - 1) -  (5N_R + 7)=  (3N_R - 8)$ degrees of freedom, which is positive when $N_R \geq 3$.

\begin{example} Suppose we have three trials. We set true parameters as 
$$
(\delta_{1, ss}, \delta_{0, ss}, \delta_{1, s\bar{s}}, \delta_{0, s\bar{s}}, \delta_{1, \bar{s}\bar{s}}, \delta_{0, \bar{s}\bar{s}}, \delta_{1, \bar{s}s}, \delta_{0, \bar{s}s})
= (0.8, 0.5, 0.7, 0.3, 0.6, 0.1, 0.5, 0.2),
$$
and display other parameters in Table B.3.1.

\begin{table}
\caption{True values of $P(Z=1 \mid R=r)$ and $P(U=u \mid R=r)$}
\\
\begin{tabular}{cccc}\hline
                                 &$r=1$ & $r=2$  & $r=3$ \\ \hline
 $P(Z=1 \mid R=r)$  & $0.4$ & $0.5$  & $0.6$ \\ 

$\pi_{ss,r}$ &$0.6$      &$0.1$                 &$0.1$     \\
$\pi_{s\bar{s},r}$ &$0.2$    &$0.6$                 &$0.1$     \\
$\pi_{\bar{s}\bar{s},r}$ &$0.1$   &$0.2$         &$0.6$    \\
$\pi_{\bar{s}s,r}$ &$0.1$     &$0.1$                   &$0.2$  \\ \hline
\end{tabular}
\end{table}

We can then get the observed distribution of $P(Z=z,S=s,Y=y\mid R=r)$ as shown in Table B.3.2.
By the homogeneity assumption, the probabilities of the observed data can be decomposed as
\begin{eqnarray}\label{eq::obs}
&&P(Z = z, S=s, Y=y\mid R=r) \nonumber \\ 
&=& \sum_{u\in O(z,s)}  P(Z=z \mid R=r) \cdot   \pi_{ur} \cdot  P(Y=y\mid Z=z, U=u),
\end{eqnarray}
for $z, s, y = 0,1$.
Similar to the proof of Proposition \ref{prop:nec}, we can obtain that the number of parameters on the right-hand side of Equation (\ref{eq::obs}) is 20. The observed distribution of $(Z,S,Y \mid R)$  form a contingency table with $21$ free frequencies on the left-hand side of Equation (\ref{eq::obs}). 
 We can use R package {\it neqslv} to solve these equations. We obtain that  $$
(\delta_{1, ss}, \delta_{0, ss}, \delta_{1, s\bar{s}}, \delta_{0, s\bar{s}}, \delta_{1, \bar{s}\bar{s}}, \delta_{0, \bar{s}\bar{s}}, \delta_{1, \bar{s}s}, \delta_{0, \bar{s}s})
= (0.8, 0.5, 0.7, 0.3, 0.6, 0.1, 0.5, 0.2),
$$  which are the same as the true parameters.
\end{example}

\begin{table} 
\caption{Observed distribution of $P(Z=z,S=s,Y=y\mid R=r)$} 
\\
\begin{tabular}{ccccccc}\hline
                    & \multicolumn{2}{c}{$R=1$} & \multicolumn{2}{c}{$R=2$}  & \multicolumn{2}{c}{$R=3$} \\
                    &  $Z=1$ & $Z=0$                 & $Z=1$ & $Z=0$                   & $Z=1$ & $Z=0$                  \\ \hline
$S=1 ,Y=1$ &0.248       &0.192                   &0.250     &0.035                   &0.090  &0.036\\
$S=1 ,Y=0$ &0.072     &0.228                  &0.100    &0.065                     &0.030    &0.084\\
$S=0 ,Y=1$ &0.044    &0.042                   &0.085    &0.100                     &0.276   &0.036\\
$S=0 ,Y=0$ &0.036     &0.138                  &0.065    &0.300                     &0.204   &0.244\\ \hline
\end{tabular}
\end{table}

\subsection*{\bf Appendix B.4: Proof of Proposition 3 in Section 4}

\begin{proposition}
\label{prop:sur:quanti}
For trial $r$, assume that the principal surrogate satisfies the causal necessity.
\begin{enumerate}
\item[(i)]
With monotonicity, we have $ACE^Y_r = ACE^S_r \times ACE_{s\bar{s},r}$.

\item[(ii)] 
Without monotonicity, suppose $ACE^S_r>0$, we have that, if $ACE_{s\bar{s},r}+ACE_{\bar{s}s,r} \geq 0$,
$$
ACE^S_r \times ACE_{s\bar{s},r} \leq ACE^Y_r \leq 
(ACE_{s\bar{s},r}+ACE_{\bar{s}s,r}) / 2 +  ACE^S_r\times  (ACE_{s\bar{s},r} - ACE_{\bar{s}s,r}) / 2,
$$
and otherwise
$$
(ACE_{s\bar{s},r}+ACE_{\bar{s}s,r}) / 2 +  ACE^S_r\times  (ACE_{s\bar{s},r} - ACE_{\bar{s}s,r}) / 2 \leq ACE^Y_r  \leq ACE^S_r \times ACE_{s\bar{s},r}.
$$

\end{enumerate}
\end{proposition}

\noindent {\it Proof:} With monotonicity, from the causal necessity $ACE_{ss,r}=ACE_{\bar{s}\bar{s},r}=0$, we have
\begin{eqnarray*}
ACE^Y_r  =  E\{Y(1)-Y(0)\mid R=r\} 
= E\{Y(1)-Y(0)\mid U=s\bar{s}, R=r\}\pi_{s\bar{s},r} 
= ACE_{s\bar{s},r} \pi_{s\bar{s},r} .
\end{eqnarray*}
Monotonicity also implies that $ACE^S_r = E\{S(1)-S(0)\mid R=r\} = \pi_{s\bar{s},r} $, and therefore $ACE^Y_r  = ACE^S_r \times ACE_{s\bar{s},r}$.

Without monotonicity, from the causal necessity, we have
$$
ACE^Y_r  = \sum_{u} E\{Y(1)-Y(0)\mid U=u, R=r\}\pi_{ur} 
=ACE_{s\bar{s}, r}\pi_{s\bar{s},r} + ACE_{\bar{s}s, r}\pi_{\bar{s}s,r}.
$$
Since $\pi_{s\bar{s},r}-\pi_{\bar{s}s,r}=E\{S(1)-S(0)\mid R=r\}=ACE^S_r$, the above equation can be rearranged as
$$
ACE^Y_r  = 
(ACE_{s\bar{s},r}+ACE_{\bar{s}s,r}) \pi_{\bar{s}s,r} +ACE^S_r \times ACE_{s\bar{s},r} .
$$
From the linear constraints $0\leq \pi_{\bar{s}s,r}, \pi_{s\bar{s},r} \leq 1, \pi_{\bar{s}s,r}+\pi_{s\bar{s},r} \leq 1$ and $\pi_{s\bar{s} ,r} - \pi_{\bar{s}s, r} = ACE^S_r$,
we obtain $\pi_{\bar{s}s,r} \in [ \max ( 0,-ACE^S_r ),(1-ACE^S_r)/2] = [0, (1-ACE^S_r)/2]$.
If $ACE_{s\bar{s},r}+ACE_{\bar{s}s,r} \geq 0$, we have
\begin{eqnarray*}
&&ACE^S_r \times ACE_{s\bar{s},r} \leq  ACE^Y_r  \\
& \leq&  (ACE_{s\bar{s},r}+ACE_{\bar{s}s,r}) \times (1-ACE^S_r)/2 +ACE^S_r \times ACE_{s\bar{s},r}\\
& =&
(ACE_{s\bar{s},r}+ACE_{\bar{s}s,r})/2 + ACE^S_r\times  (ACE_{s\bar{s},r} - ACE_{\bar{s}s,r}) / 2.
\end{eqnarray*}
Otherwise, we have 
$$
(ACE_{s\bar{s},r}+ACE_{\bar{s}s,r}) / 2 +  ACE^S_r\times  (ACE_{s\bar{s},r} - ACE_{\bar{s}s,r}) / 2 \leq ACE^Y_r  \leq ACE^S_r \times ACE_{s\bar{s},r}.
$$

\section*{\bf Appendix C: EM algorithms and the Gibbs Samplers}
\setcounter{equation}{0}
\renewcommand {\theequation} {C.\arabic{equation}}
Let $N_{total}$ denote the total sample size.
Let $n_{zuyr}$ be the frequencies of the complete data, i.e., the frequencies of patients with
$
(Z=z, U=u, Y=y, R=r)$,
for $z=0,1$; $r=1,\ldots,N_R$; $y=0,1$; and $u=ss, s\bar{s}, \bar{s}\bar{s}, \bar{s}s$.
Let $N_{zsyr}$ be the frequencies of the observed data, i.e., the frequencies of patients with
$
(Z=z, S=s, Y=y, R=r)$,
for $z=0,1$; $r=1,\ldots,N_R$; $y=0,1$; and $s=0,1$.

\subsection*{\bf Appendix C.1: EM algorithm with monotonicity}
We can write the complete-data likelihood as
\begin{align}
L_c(\boldsymbol{\xi})  \propto & \prod^{N_R}_{r=1}\bigg(p_{r}^{n_{+++r}} \cdot \alpha_r^{n_{1++r}} (1-\alpha_r)^{n_{0++r}} \cdot \prod_{u=ss,s\bar{s},\bar{s}\bar{s}}\pi_{ur}^{n_{+u+r}} \bigg)    \nonumber \\
&
\cdot \prod_{z,u} \delta_{zu}^{n_{zu1+}}(1-\delta_{zu})^{n_{zu0+}}\cdot \prod_{i=1}^{N_{total}} I\{U_i \in O(Z_i,S_i)\}, \label{eq::comlikmono}
\end{align}
where $I\{\cdot\}$ is the indicator function.
Let $\boldsymbol{\xi}^{(k)}$ be the estimate of $\boldsymbol{\xi}$ after the $k$-th iteration.

In the E-step, we need to calculate $n_{zuyr}^{(k+1)}=E(n_{zuyr}\mid N_{zsyr}, \boldsymbol{\xi}^{(k)})$ given the observed data and $\boldsymbol{\xi}^{(k)}$.
Under monotonicity, it is relatively straightforward to obtain that $n_{0,ss,y,r}^{(k+1)}= N_{01yr}$ and $n_{1,\bar{s}\bar{s},y,r}^{(k+1)}= N_{10yr}$.
For other frequencies, we have
\begin{eqnarray*}
n_{1,ss,1,r}^{(k+1)}=&N_{111r}\frac{\pi_{ss,r}^{(k)} \delta_{1,ss}^{(k)}}{\pi_{ss,r}^{(k)} \delta_{1,ss}^{(k)}+\pi_{s\bar{s},r}^{(k)} \delta_{1,s\bar{s}}^{(k)}},& n_{1,s\bar{s},1,r}^{(k+1)}=N_{111r}-n_{1,ss,1,r}^{(k+1)} ,\\
n_{1,ss,0,r}^{(k+1)}=&N_{110r}\frac{\pi_{ss,r}^{(k)}(1- \delta_{1,ss}^{(k)})}{\pi_{ss,r}^{(k)} (1-\delta_{1,ss}^{(k)})+\pi_{s\bar{s},r}^{(k)} (1-\delta_{1,s\bar{s}}^{(k)})},& n_{1,s\bar{s},0,r}^{(k+1)}=N_{110r}-n_{1,ss,0,r}^{(k+1)},\\
n_{0,\bar{s}\bar{s},1,r}^{(k+1)}=&N_{001r}\frac{\pi_{\bar{s}\bar{s},r}^{(k)} \delta_{1,\bar{s}\bar{s}}^{(k)}}{\pi_{\bar{s}\bar{s},r}^{(k)} \delta_{1,\bar{s}\bar{s}}^{(k)}+\pi_{s\bar{s},r}^{(k)} \delta_{1,s\bar{s}}^{(k)}}, &n_{0,\bar{s}s,1,r}^{(k+1)}=N_{001r}-n_{0,\bar{s}\bar{s},1,r}^{(k+1)} ,\\
n_{0,\bar{s}\bar{s},0,r}^{(k+1)}=&N_{000r}\frac{\pi_{\bar{s}\bar{s},r}^{(k)}(1- \delta_{1,\bar{s}\bar{s}}^{(k)})}{\pi_{\bar{s}\bar{s},r}^{(k)} (1-\delta_{1,\bar{s}\bar{s}}^{(k)})+\pi_{s\bar{s},r}^{(k)} (1-\delta_{1,s\bar{s}}^{(k)})}, & n_{0,\bar{s}s,0,r}^{(k+1)}=N_{000r}-n_{0,\bar{s}\bar{s},0,r}^{(k+1)}.
\end{eqnarray*}

In the M-step, we update the parameters as follows:
\begin{eqnarray*}
p_r^{(k+1)} = \frac{n_{+++r}^{(k+1)}}{n_{++++}^{(k+1)}}, 
\quad 
\alpha_r^{(k+1)}=\frac{n_{1++r}^{(k+1)}}{n_{+++r}^{(k+1)}}, 
\quad
\pi_{ur}^{(k+1)}=\frac{n_{+u+r}^{(k+1)}}{n_{+++r}^{(k+1)}}, 
\quad
\delta_{zu}^{(k+1)}=\frac{n_{zu1+}^{(k+1)}}{n_{zu++}^{(k+1)}}.
\end{eqnarray*}

\subsection*{\bf Appendix C.2: EM algorithm without monotonicity}
We can write the complete-data likelihood as
\begin{align}
L_c(\boldsymbol{\xi})  \propto & \prod^{N_R}_{r=1}\bigg(p_{r}^{n_{+++r}} \cdot \alpha_r^{n_{1++r}} (1-\alpha_r)^{n_{0++r}} \cdot \prod_{u=ss,s\bar{s},\bar{s}\bar{s},\bar{s}s}\pi_{ur}^{n_{+u+r}} \bigg)  \nonumber \\
&
\cdot  \prod_{z,u} \delta_{zu}^{n_{zu1+}}(1-\delta_{zu})^{n_{zu0+}}\cdot \prod_{i=1}^{N_{total}} I\{U_i \in O(Z_i,S_i)\}, \label{eq::comliknmono}
\end{align}

In the E-step, we need to calculate $n_{zuyr}^{(k+1)}=E(n_{zuyr}\mid N_{zsyr}, \boldsymbol{\xi}^{(k)})$ given the observed data and $\boldsymbol{\xi}^{(k)}$ as follows:
\begin{eqnarray*}
&n_{1,ss,1,r}^{(k+1)}=N_{111r}\frac{\pi_{ss,r}^{(k)} \delta_{1,ss}^{(k)}}{\pi_{ss,r}^{(k)} \delta_{,1ss}^{(k)}+\pi_{s\bar{s},r}^{(k)} \delta_{1,s\bar{s}}^{(k)}}, &n_{1,s\bar{s},1,r}^{(k+1)}=N_{111r}-n_{1,ss,1,r}^{(k+1)} ,\\
&n_{1,ss,0,r}^{(k+1)}=N_{110r}\frac{\pi_{ss,r}^{(k)}(1- \delta_{1,ss}^{(k)})}{\pi_{ss,r}^{(k)} (1-\delta_{1,ss}^{(k)})+\pi_{s\bar{s},r}^{(k)} (1-\delta_{1,s\bar{s}}^{(k)})}, &n_{1,s\bar{s},0,r}^{(k+1)}=N_{110r}-n_{1,ss,0,r}^{(k+1)} ,\\
&n_{1,\bar{s}\bar{s},1,r}^{(k+1)}=N_{101r}\frac{\pi_{\bar{s}\bar{s},r}^{(k)} \delta_{1,\bar{s}\bar{s}}^{(k)}}{\pi_{\bar{s}\bar{s},r}^{(k)} \delta_{1,\bar{s}\bar{s}}^{(k)}+\pi_{\bar{s}s,r}^{(k)} \delta_{1,\bar{s}s}^{(k)}}, &n_{1,\bar{s}s,1,r}^{(k+1)}=N_{101r}-n_{1,\bar{s}\bar{s},1,r}^{(k+1)} ,\\
&n_{1,\bar{s}\bar{s},0,r}^{(k+1)}=N_{100r}\frac{\pi_{\bar{s}\bar{s},r}^{(k)}(1- \delta_{1,\bar{s}\bar{s}}^{(k)})}{\pi_{\bar{s}\bar{s},r}^{(k)} (1-\delta_{1,\bar{s}\bar{s}}^{(k)})+\pi_{\bar{s}s,r}^{(k)} (1-\delta_{1,\bar{s}s}^{(k)})}, &n_{1,\bar{s}s,0,r}^{(k+1)}=N_{100r}-n_{1,\bar{s}\bar{s},0,r}^{(k+1)} ,\\
&n_{0,ss,1,r}^{(k+1)}=N_{011r}\frac{\pi_{ss,r}^{(k)} \delta_{1,ss}^{(k)}}{\pi_{ss,r}^{(k)} \delta_{1,ss}^{(k)}+\pi_{\bar{s}s,r}^{(k)} \delta_{1,\bar{s}s}^{(k)}}, &n_{0,\bar{s}s,1,r}^{(k+1)}=N_{011r}-n_{0,ss,1,r}^{(k+1)} ,\\
&n_{0,ss,0,r}^{(k+1)}=N_{010r}\frac{\pi_{ss,r}^{(k)}(1- \delta_{1,ss}^{(k)})}{\pi_{ss,r}^{(k)} (1-\delta_{1,ss}^{(k)})+\pi_{\bar{s}s,r}^{(k)} (1-\delta_{1,\bar{s}s}^{(k)})}, &n_{0,\bar{s}s,0,r}^{(k+1)}=N_{010r}-n_{0,ss,0,r}^{(k+1)} ,\\
&n_{0,\bar{s}\bar{s},1,r}^{(k+1)}=N_{001r}\frac{\pi_{\bar{s}\bar{s},r}^{(k)} \delta_{1,\bar{s}\bar{s}}^{(k)}}{\pi_{\bar{s}\bar{s},r}^{(k)} \delta_{1,\bar{s}\bar{s}}^{(k)}+\pi_{s\bar{s},r}^{(k)} \delta_{1,s\bar{s}}^{(k)}}, &n_{0,\bar{s}s,1,r}^{(k+1)}=N_{001r}-n_{0,\bar{s}\bar{s},1,r}^{(k+1)} ,\\
&n_{0,\bar{s}\bar{s},0,r}^{(k+1)}=N_{000r}\frac{\pi_{\bar{s}\bar{s},r}^{(k)}(1- \delta_{1,\bar{s}\bar{s}}^{(k)})}{\pi_{\bar{s}\bar{s},r}^{(k)} (1-\delta_{1,\bar{s}\bar{s}}^{(k)})+\pi_{s\bar{s},r}^{(k)}(1-\delta_{1,s\bar{s}}^{(k)})},
&n_{0,\bar{s}s,0,r}^{(k+1)}=N_{000r}-n_{0,\bar{s}\bar{s},0,r}^{(k+1)}.
\end{eqnarray*}

In the M-step, we can update the parameters as follows:
\begin{eqnarray*}
p_r^{(k+1)} = \frac{n_{+++r}^{(k+1)}}{n_{++++}^{(k+1)}},
\quad
 \alpha_r^{(k+1)}=\frac{n_{1++r}^{(k+1)}}{n_{+++r}^{(k+1)}}, 
 \quad
 \pi_{ur}^{(k+1)}=\frac{n_{+u+r}^{(k+1)}}{n_{+++r}^{(k+1)}}, 
 \quad
 \delta_{zu}^{(k+1)}=\frac{n_{zu1+}^{(k+1)}}{n_{zu++}^{(k+1)}}.
\end{eqnarray*}

\subsection*{\bf Appendix C.3: Gibbs Sampler with monotonicity}
We treat $U$ as the missing data, and the Gibbs sampler iterates between the following imputation and posterior steps.
In the imputation step, we draw $U$ given the observed data and all the parameters.
When $(Z_i=1,S_i=0)$, we impute $U_i=\bar{s}\bar{s}$; when $(Z_i=0,S_i=1)$, we impute $U_i=ss$; when $(Z_i=1,S_i=1)$, we impute $U_i$ as:
\begin{eqnarray*}
&&P(U_i=s\bar{s}\mid Z_i=1,S_i=1,R_i,Y_i,\bm{\xi})  = \frac{\delta_{1,s\bar{s}}^{Y_i}(1-\delta_{1,s\bar{s}})^{1-Y_i}\pi_{s\bar{s},R_i}}{\sum_{u=ss,s\bar{s}}\delta_{1,u}^{Y_i}(1-\delta_{1,u})^{1-Y_i}\pi_{u,R_i}},\\
&&P(U_i=ss \mid Z_i=1,S_i=1,R_i,Y_i,\bm{\xi}) = 1- P(U_i=s\bar{s} \mid Z_i=1,S_i=1,R_i,Y_i,\bm{\xi});
\end{eqnarray*}
when $(Z_i=0,S_i=0)$, we impute $U_i$ as:
\begin{eqnarray*}
&&P(U_i=s\bar{s} \mid Z_i=0,S_i=0,R_i,Y_i,\bm{\xi})  = \frac{\delta_{0,s\bar{s}}^{Y_i}(1-\delta_{0,s\bar{s}})^{1-Y_i}\pi_{s\bar{s},R_i}}{\sum_{u=\bar{s}\bar{s},s\bar{s}}\delta_{0,u}^{Y_i}(1-\delta_{0,u})^{1-Y_i}\pi_{u,R_i}},\\
&&P(U_i=\bar{s}\bar{s} \mid Z_i=0,S_i=0,R_i,Y_i,\bm{\xi}) = 1- P(U_i=s\bar{s} \mid Z_i=0,S_i=0,R_i,Y_i,\bm{\xi}).
\end{eqnarray*}
After imputing $U$, we can compute the frequencies of the complete data.

In the posterior step, we draw all the parameters given the complete data. Specifically, we draw the parameters from the following conditional distributions:
\begin{eqnarray*}
(p_1,p_2,\cdots,p_{N_R} )  \mid \cdot &\sim& \text{Dirichlet}(n_{+++1}+1, \cdots, n_{+++N_R}+1),\\
\alpha_r \mid \cdot &\sim& \text{Beta}(n_{1++r}+1,n_{0++r}+1), \forall r, \\
(  \pi_{ss,r},\pi_{s\bar{s},r},\pi_{\bar{s}\bar{s},r} )  \mid \cdot &\sim& \text{Dirichlet}(n_{+ss+r}+1,n_{+s\bar{s}+r}+1,n_{+\bar{s}\bar{s}+r}+1),\\
\delta_{zu} \mid \cdot &\sim& \text{Beta}(n_{zu1+}+1,n_{zu0+}+1), \forall u, z.
\end{eqnarray*}

\subsection*{\bf Appendix C.4: Gibbs Sampler without monotonicity}

We treat $U$ as the missing data, and the Gibbs sampler iterates between the following imputation and posterior steps.
In the imputation step, we draw $U$ given the observed data and all the parameters. When $(Z_i=1,S_i=1)$, we impute $U_i$ as:
\begin{eqnarray*}
&&P(U_i=s\bar{s} \mid Z_i=1,S_i=1,R_i,Y_i,\bm{\xi})  = \frac{\delta_{1,s\bar{s}}^{Y_i}(1-\delta_{1,s\bar{s}})^{1-Y_i}\pi_{s\bar{s},R_i}}{\sum_{u=ss,s\bar{s}}\delta_{1,u}^{Y_i}(1-\delta_{1,u})^{1-Y_i}\pi_{u,R_i}},\\
&&P(U_i=ss \mid Z_i=1,S_i=1,R_i,Y_i,\bm{\xi}) = 1- P(U_i=s\bar{s} \mid Z_i=1,S_i=1,R_i,Y_i,\bm{\xi}); 
\end{eqnarray*}
when $(Z_i=0,S_i=0)$, we impute $U_i$ as:
\begin{eqnarray*}
&&P(U_i=s\bar{s} \mid Z_i=0,S_i=0,R_i,Y_i,\bm{\xi})  = \frac{\delta_{0,s\bar{s}}^{Y_i}(1-\delta_{0,s\bar{s}})^{1-Y_i}\pi_{s\bar{s},R_i}}{\sum_{u=\bar{s}\bar{s},s\bar{s}}\delta_{0,u}^{Y_i}(1-\delta_{0,u})^{1-Y_i}\pi_{u,R_i}},\\
&&P(U_i=\bar{s}\bar{s} \mid Z_i=0,S_i=0,R_i,Y_i,\bm{\xi}) = 1- P(U_i=s\bar{s} \mid Z_i=0,S_i=0,R_i,Y_i,\bm{\xi});
\end{eqnarray*}
when $(Z_i=1,S_i=0)$, we impute $U_i$ as:
\begin{eqnarray*}
&&P(U_i=\bar{s}s \mid Z_i=1,S_i=0,R_i,Y_i,\bm{\xi})  = \frac{\delta_{1,\bar{s}s}^{Y_i}(1-\delta_{1,\bar{s}s})^{1-Y_i}\pi_{\bar{s}s,R_i}}{\sum_{u=\bar{s}s,\bar{s}\bar{s}}\delta_{1,u}^{Y_i}(1-\delta_{1,u})^{1-Y_i}\pi_{u,R_i}},\\
&&P(U_i=\bar{s}\bar{s} \mid Z_i=1,S_i=0,R_i,Y_i,\bm{\xi}) = 1- P(U_i=\bar{s}s \mid Z_i=1,S_i=0,R_i,Y_i,\bm{\xi});
\end{eqnarray*}
when $(Z_i=0,S_i=1)$, we impute $U_i$ as:
\begin{eqnarray*}
&&P(U_i=\bar{s}s \mid Z_i=0,S_i=1,R_i,Y_i,\bm{\xi})  = \frac{\delta_{0,\bar{s}s}^{Y_i}(1-\delta_{0,\bar{s}s})^{1-Y_i}\pi_{\bar{s}s,R_i}}{\sum_{u=\bar{s}s,ss}\delta_{0,u}^{Y_i}(1-\delta_{0,u})^{1-Y_i}\pi_{u,R_i}},\\
&&P(U_i=ss \mid Z_i=0,S_i=1,R_i,Y_i,\bm{\xi}) = 1- P(U_i=\bar{s}s \mid Z_i=0,S_i=1,R_i,Y_i,\bm{\xi}).
\end{eqnarray*}
After imputing $U$, we can compute the frequencies of the complete data.

In the posterior step, we draw all the parameters given the complete data. Specifically, we draw the parameters from the following conditional distributions:
\begin{eqnarray*}
 ( p_1,p_2,\cdots,p_{N_R} )   \mid \cdot &\sim& \text{Dirichlet}(n_{+++1}+1, \cdots, n_{+++N_R}+1),\\
\alpha_r \mid \cdot &\sim& \text{Beta}(n_{1++r}+1,n_{0++r}+1), \forall r,\\
(  \pi_{ss,r},\pi_{s\bar{s},r},\pi_{\bar{s}\bar{s},r},\pi_{\bar{s}s,r} )  \mid \cdot &\sim& \text{Dirichlet}(n_{+ss+r}+1,n_{+s\bar{s}+r}+1,n_{+\bar{s}\bar{s}+r}+1,n_{+\bar{s}s+r}+1),\\
\delta_{zu} \mid \cdot &\sim& \text{Beta}(n_{zu1+}+1,n_{zu0+}+1), \forall z, u.
\end{eqnarray*}

\subsection*{\bf Appendix C.5: Gibbs Sampler for the Bayesian hierarchical model in Section 7}

Define $\eta_{zur} = \text{logit} (\delta_{zur}), \bm{\eta} = \{ \eta_{zur}: z=0,1; u=ss, s\bar{s}, \bar{s}s, \bar{s}\bar{s}   \}$, and $\boldsymbol{\theta} = (\bm{\xi}, \bm{\eta})$. The complete data likelihood is
\begin{eqnarray*}
&&L_c(\bm{\theta}) \\
& \propto & \prod^{N_R}_{r=1}\left\{p_{r}^{n_{+++r}} \cdot \alpha_r^{n_{1++r}} (1-\alpha_r)^{n_{0++r}}
       \prod_{u=ss,s\bar{s},\bar{s}\bar{s},\bar{s}s}\pi_{ur}^{n_{+u+r}} \right\}\\
&&\cdot  \prod_{z,u,r} \left [\{\text{logit}^{-1}(\eta_{zur})\}^{n_{zu1r}}\{1-\text{logit}^{-1}(\eta_{zur})\}^{n_{zu0r}}
 (\sigma^2)^{-1/2}
 \exp(-\frac{(\eta_{zur}-\mu_{zu})^2}{2\sigma^2})  \right]  \\
&&\cdot  \prod_{i=1}^{N_{total}} I\{U_i \in O(Z_i,S_i)\}.
\end{eqnarray*}

We treat $U$ as the missing data, and the Gibbs sampler iterates between the following imputation and posterior steps.
In the imputation step, we draw $U$ given the observed data and all the parameters.
When $(Z_i=1,S_i=1)$, we impute $U_i$ as:
\begin{eqnarray*}
&&P(U_i=s\bar{s} \mid Z_i=1,S_i=1,R_i,Y_i,\bm{\xi})  = \frac{\delta_{1,s\bar{s},R_i}^{Y_i}(1-\delta_{1,s\bar{s},R_i})^{1-Y_i}\pi_{s\bar{s},R_i}}{\sum_{u=ss,s\bar{s}}\delta_{1,u,R_i}^{Y_i}(1-\delta_{1,u,R_i})^{1-Y_i}\pi_{u,R_i}},\\
&&P(U_i=ss \mid Z_i=1,S_i=1,R_i,Y_i,\bm{\xi}) = 1- P(U_i=s\bar{s} \mid Z_i=1,S_i=1,R_i,Y_i,\bm{\xi});
\end{eqnarray*}
when $(Z_i=0,S_i=0)$, we impute $U_i$ as:
\begin{eqnarray*}
&&P(U_i=s\bar{s} \mid Z_i=0,S_i=0,R_i,Y_i,\bm{\xi})  = \frac{\delta_{0,s\bar{s},R_i}^{Y_i}(1-\delta_{0,s\bar{s},R_i})^{1-Y_i}\pi_{s\bar{s},R_i}}{\sum_{u=\bar{s}\bar{s},s\bar{s}}\delta_{0,u,R_i}^{Y_i}(1-\delta_{0,u,R_i})^{1-Y_i}\pi_{u,R_i}},\\
&&P(U_i=\bar{s}\bar{s} \mid Z_i=0,S_i=0,R_i,Y_i,\bm{\xi}) = 1- P(U_i=s\bar{s} \mid Z_i=0,S_i=0,R_i,Y_i,\bm{\xi});
\end{eqnarray*}
when $(Z_i=1,S_i=0)$, we impute $U_i$ as:
\begin{eqnarray*}
&&P(U_i=\bar{s}s \mid Z_i=1,S_i=0,R_i,Y_i,\bm{\xi})  = \frac{\delta_{1,\bar{s}s,R_i}^{Y_i}(1-\delta_{1,\bar{s}s,R_i})^{1-Y_i}\pi_{\bar{s}s,R_i}}{\sum_{u=\bar{s}s,\bar{s}\bar{s}}\delta_{1,u,R_i}^{Y_i}(1-\delta_{1,u,R_i})^{1-Y_i}\pi_{u,R_i}},\\
&&P(U_i=\bar{s}\bar{s} \mid Z_i=1,S_i=0,R_i,Y_i,\bm{\xi}) = 1- P(U_i=\bar{s}s \mid Z_i=1,S_i=0,R_i,Y_i,\bm{\xi});
\end{eqnarray*}
when $(Z_i=0,S_i=1)$, we impute $U_i$ as:
\begin{eqnarray*}
&&P(U_i=\bar{s}s \mid Z_i=0,S_i=1,R_i,Y_i,\bm{\xi})  = \frac{\delta_{0,\bar{s}s,R_i}^{Y_i}(1-\delta_{0,\bar{s}s,R_i})^{1-Y_i}\pi_{\bar{s}s,R_i}}{\sum_{u=\bar{s}s,ss}\delta_{0,u,R_i}^{Y_i}(1-\delta_{0,u,R_i})^{1-Y_i}\pi_{u,R_i}},\\
&&P(U_i=ss \mid Z_i=0,S_i=1,R_i,Y_i,\bm{\xi}) = 1- P(U_i=\bar{s}s \mid Z_i=0,S_i=1,R_i,Y_i,\bm{\xi}).
\end{eqnarray*}
After imputing $U$, we can compute the frequencies of the complete data.

In the posterior step, we draw all the parameters given the complete data. Specifically, we draw the parameters from the following conditional distributions:
\begin{eqnarray*}
( p_1,p_2,\cdots,p_{N_R} )  \mid \cdot &\sim &\text{Dirichlet}(n_{+++1}+1, \cdots, n_{+++N_R}+1),\\
\alpha_r \mid \cdot &\sim& \text{Beta}(n_{1++r}+1,n_{0++r}+1), \forall r,\\
( \pi_{ss,r},\pi_{s\bar{s},r},\pi_{\bar{s}\bar{s},r},\pi_{\bar{s}s,r})  \mid \cdot
&\sim& \text{Dirichlet}(n_{+ss+r}+1,n_{+s\bar{s}+r}+1,n_{+\bar{s}\bar{s}+r}+1,n_{+\bar{s}s+r}+1),\\
\mu_{zu} \mid \cdot &\propto& \exp\left\{\frac{\sum_{r=1}^{N_R}(\eta_{zur}-\mu_{zu})^2}{2\sigma^2}\right\}\cdot I(\mu_{zu} \in [-5,5]), \forall z, u, \\
p(\eta_{zur} \mid \cdot) &\propto& \{\text{expit}(\eta_{zur})\}^{n_{zu1r}}\{1-\text{expit}(\eta_{zur})\}^{n_{zu0r}} \\
&& \cdot  \exp\left\{-\frac{(\eta_{zur}-\mu_{zu})^2}{2\sigma^2}\right\}  \\
& \propto& \frac{e^{ -n_{zu0r} \cdot \eta_{zur}-\frac{(\eta_{zur}-\mu_{zu})^2}{2\sigma^2} }}{(1+e^{-\eta_{zur}})^{n_{zu+r}}}, \forall z, u,r.
\end{eqnarray*}

The first three conditional distributions are standard. The fourth one is a truncated Normal distribution, which can be generated by applying the inverse of its cumulative distribution function to a Uniform$(0,1)$ random variable.
The posterior distribution of $\eta_{zur}$ is not standard, but fortunately its conditional density is log concave since
\begin{eqnarray*}
\frac{\partial^2 \log\{p(\eta_{zur} \mid \cdot)\}}{\partial \eta_{zur}^2} = -\frac{1}{\sigma^2}-\frac{n_{zu+r} e^{-\eta_{zur}}}{(1+e^{-\eta_{zur}})^2}<0.
\end{eqnarray*}
Therefore its posterior distribution is unimodal, and we can use the Metropolized Independence Sampler ~\citep{liu2001generalized} to sample $\eta_{zur}$ using a Normal proposal. We choose the proposal distribution as
$
\eta \sim N(\widehat{\eta},\widehat{\psi}),
$
where $\widehat{\eta}$ is the mode of the posterior distribution and $\widehat{\psi}$ is the inverse of the Fisher information at the mode. Due to the log concave density, we can simply use the Newton-Raphson iteration to find the mode $\widehat{\eta}$.

\section*{\bf Appendix D: More details about the simulation studies}

\setcounter{equation}{0}
\renewcommand {\theequation} {D.\arabic{equation}}

For the correctly specified model, Figures \ref{fig:sim1:mon:HL} and \ref{fig:sim1:non:HL} show the results for $ACE_{s\bar{s}}$, and Figures \ref{fig:sim1:mon:LL} and \ref{fig:sim1:non:LL} show the results for $ACE_{\bar{s}\bar{s}}$, with the biases and  RMSEs of the MLEs and the coverage proportions of the posterior credible intervals. For the misspecified models, Figures \ref{fig:sim2:mon:HL} and \ref{fig:sim2:non:HL} show the results for $ACE_{s\bar{s}}$, Figures \ref{fig:sim2:mon:LL} and \ref{fig:sim2:non:LL} show the results for $ACE_{\bar{s}\bar{s}}$, and Figures \ref{fig:sim2:mon:LH} and \ref{fig:sim2:non:LH} show the results for $ACE_{\bar{s}s}$, with the biases and RMSEs of the MLEs and the coverage proportions of the posterior credible intervals.
The reason why the RMSEs with more trials may be larger for the same strata is that certain settings of parameters may make the sample sizes smaller for the strata with more trials.

\setcounter{figure}{4}
\begin{figure}
  \centering
\subfigure[With monotonicity]{
    \label{fig:2b} 
\includegraphics[width=\textwidth]{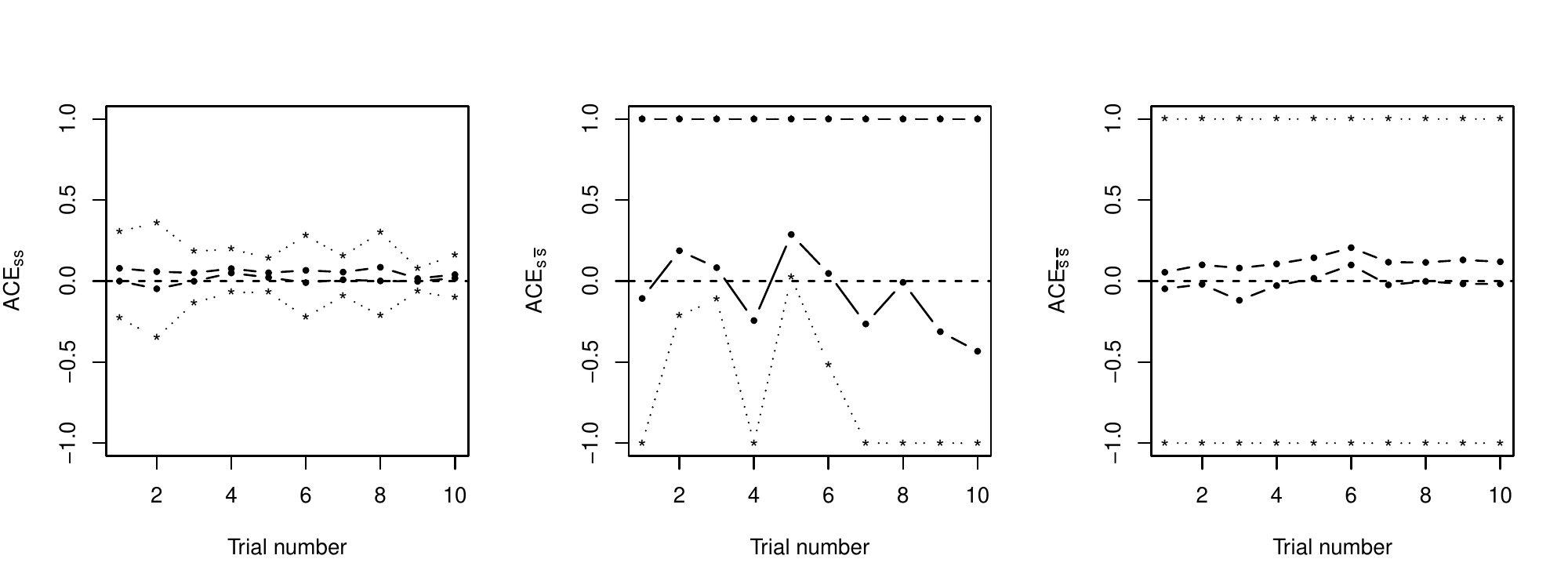}}
  \subfigure[Without monotonicity]{
    \label{fig:2a} 
\includegraphics[width=\textwidth]{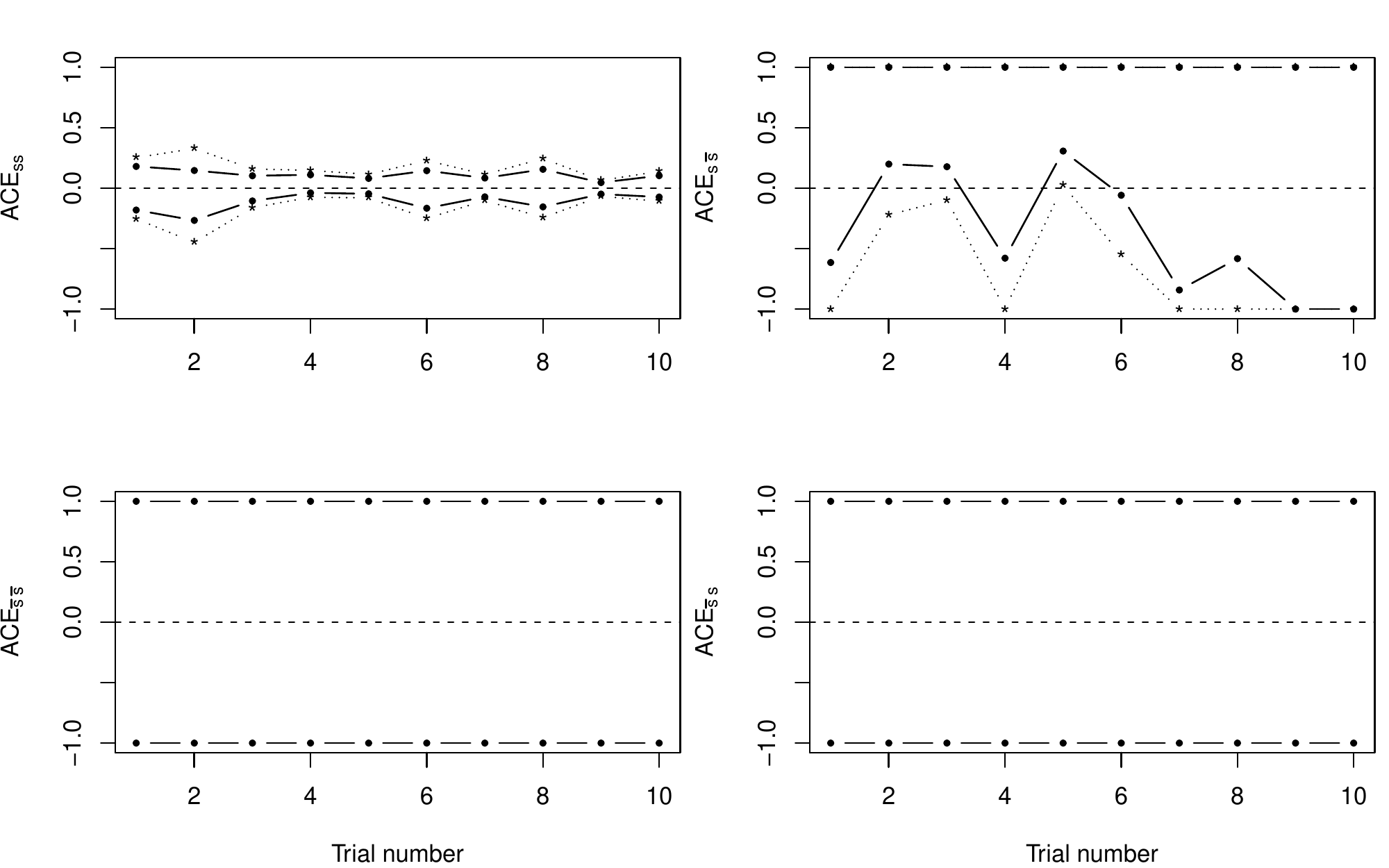}}
  \caption{Bounds for the PSACEs with and without monotonicity. The solid lines are the bounds and the dotted lines are the confidence intervals for the bounds.}
  \label{fig:2} 
\end{figure}

\begin{figure}
  \centering
   \subfigure[Correctly specified model with monotonicity]{
   \label{fig:sim1:mon:HL}
\scalebox{1}[0.95]{\includegraphics[width=\textwidth]{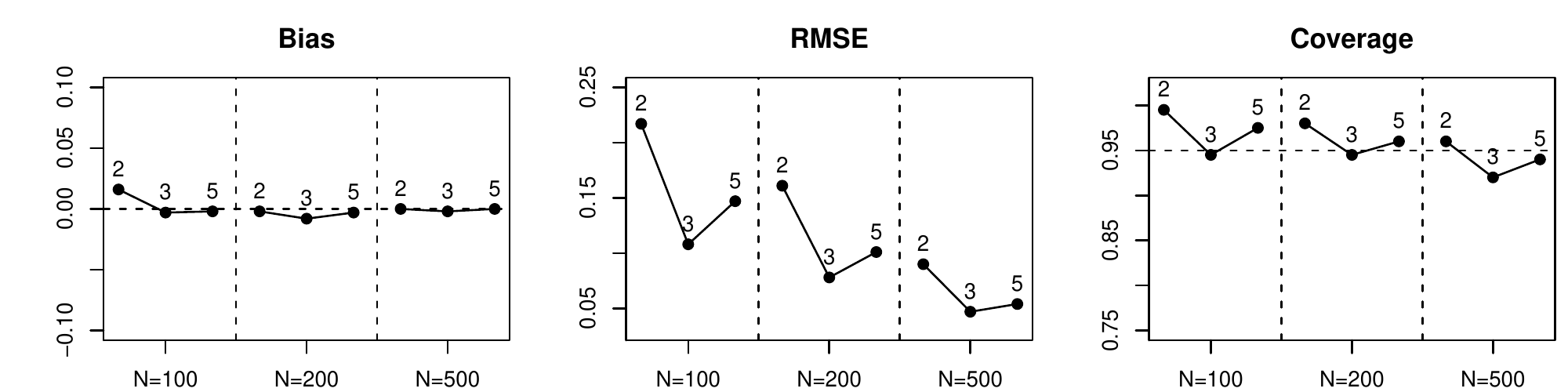}}} \subfigure[Correctly specified model without monotonicity]{
 \label{fig:sim1:non:HL}
\scalebox{1}[0.95]{\includegraphics[width=\textwidth]{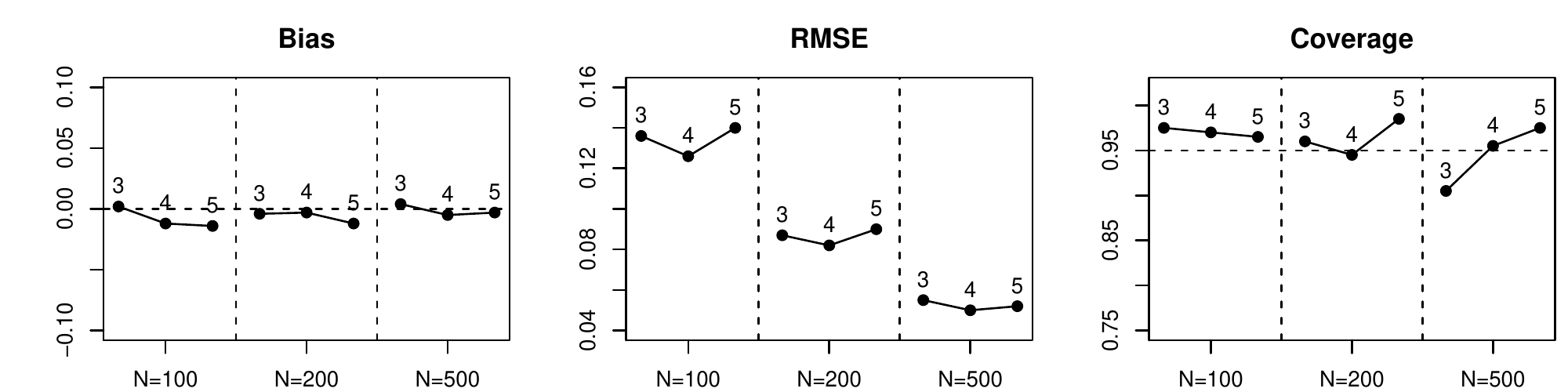}}}
  \subfigure[Misspecified model with monotonicity]{
  \label{fig:sim2:mon:HL}
\scalebox{1}[0.95]{\includegraphics[width=\textwidth]{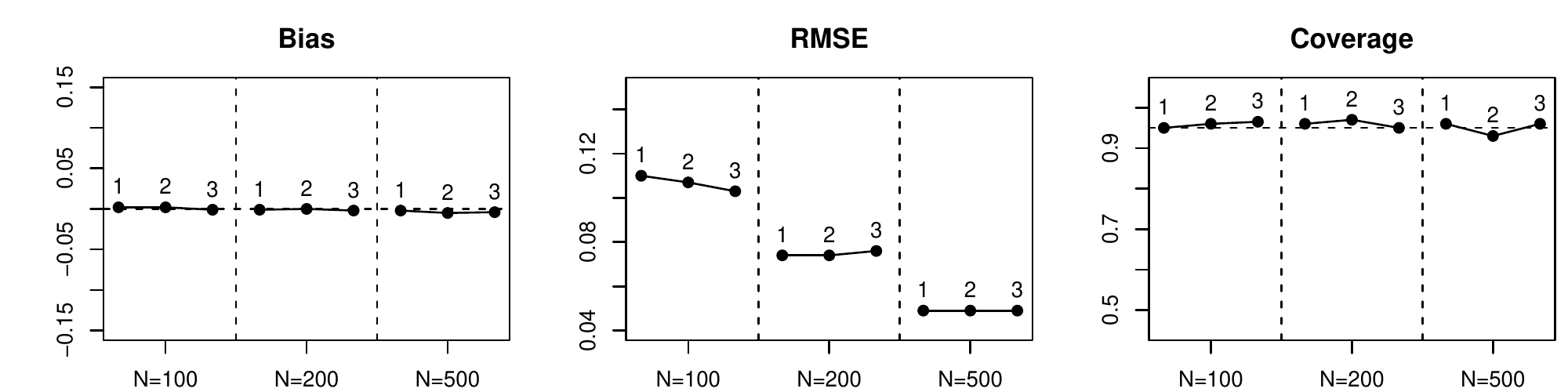}}} \subfigure[Misspecified model without monotonicity]{
  \label{fig:sim2:non:HL}
\scalebox{1}[0.95]{\includegraphics[width=\textwidth]{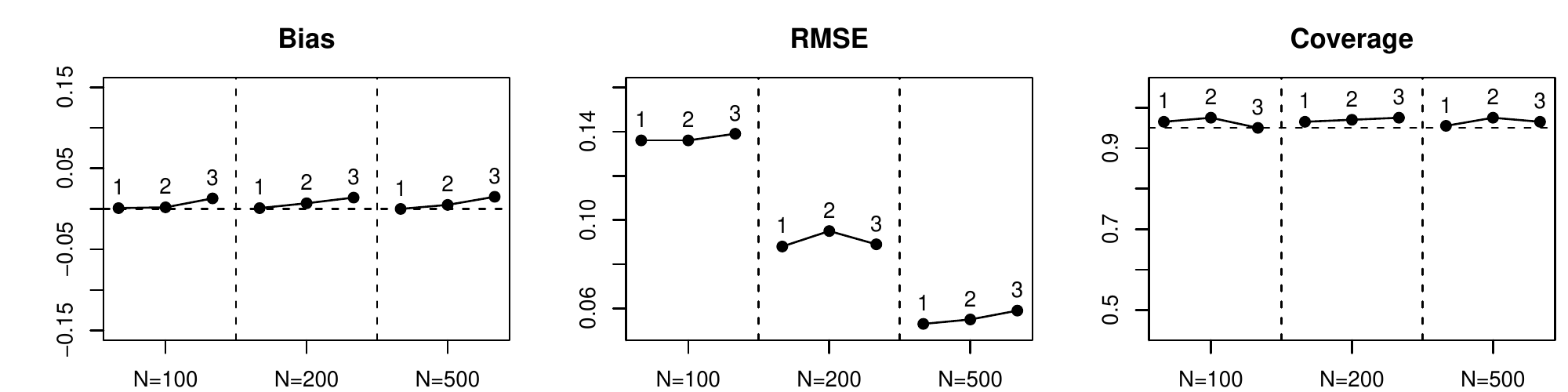}}}
  \caption{Simulations for homogeneity (a) and (b) and heterogeneity (c) and (d). Each subgraph presents the bias, RMSE, coverage proportions of the $95\%$ credible intervals of  $ACE_{s\bar{s}}$. Nine combinations of $N_R$ and $N$ are shown in (a) and (b); and nine combinations of $N$ and $d$ (``1'' for $.01$, ``2'' for $.025$, ``3'' for $.05$) are shown in (c) and (d).}
\end{figure}

\begin{figure}
  \centering
   \subfigure[Correctly specified model with monotonicity]{
   \label{fig:sim1:mon:LL}
\scalebox{1}[0.95]{\includegraphics[width=\textwidth]{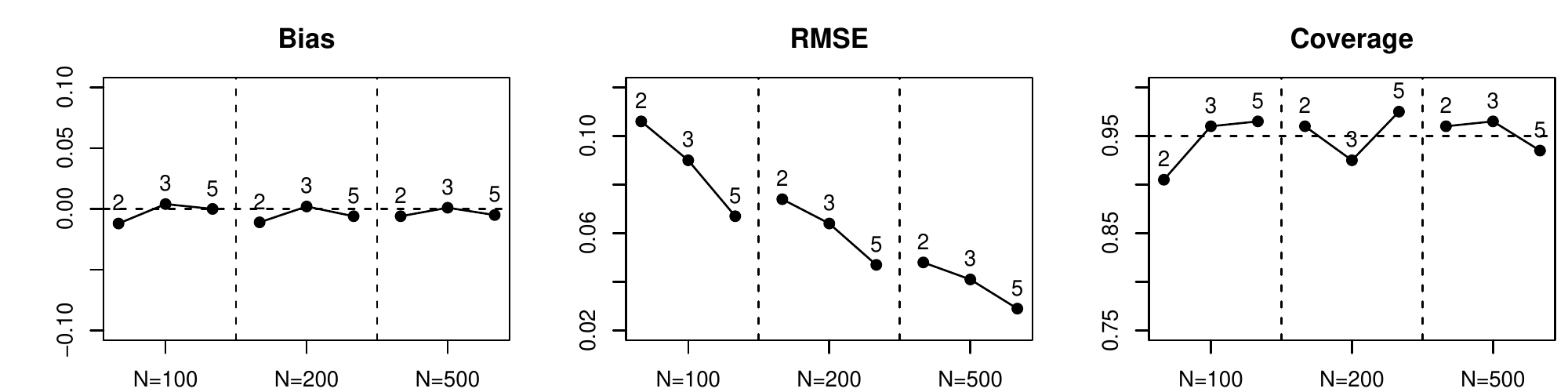}}} \subfigure[Correctly specified model without monotonicity]{
\label{fig:sim1:non:LL}
\scalebox{1}[0.95]{\includegraphics[width=\textwidth]{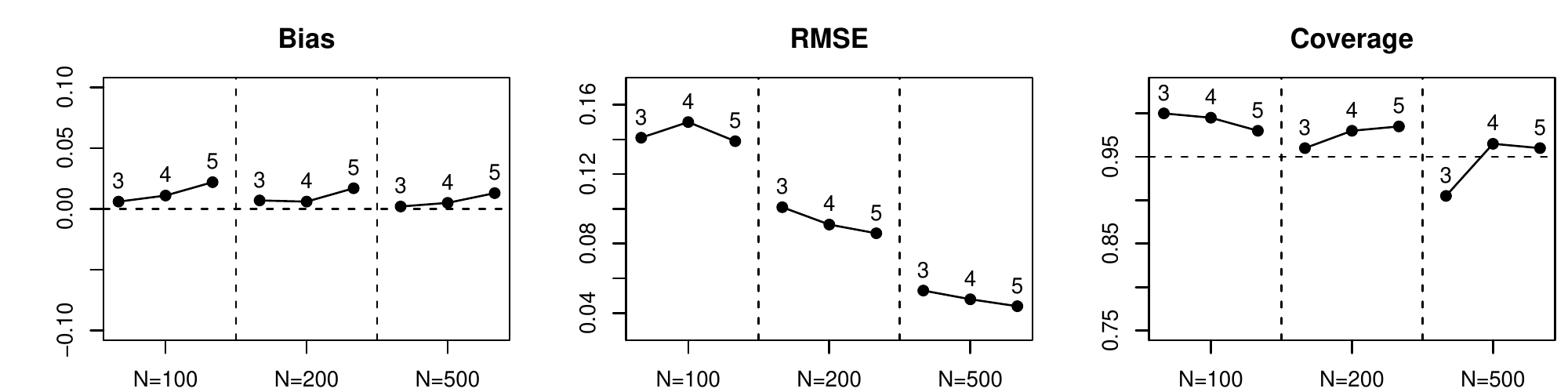}}}
  \subfigure[Misspecified model with monotonicity]{
  \label{fig:sim2:mon:LL}
\scalebox{1}[0.95]{\includegraphics[width=\textwidth]{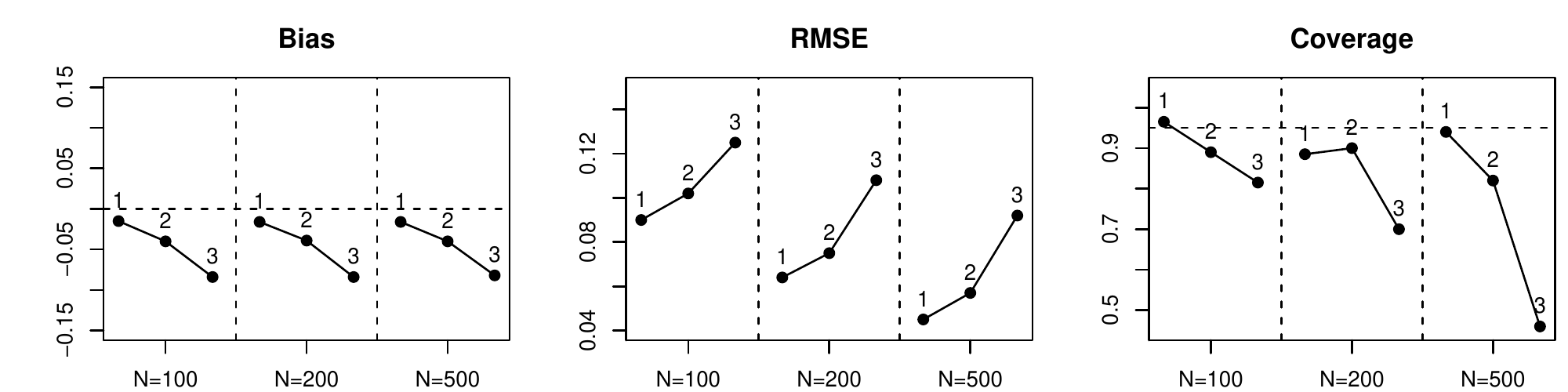}}} \subfigure[Misspecified model without monotonicity]{
\label{fig:sim2:non:LL}
\scalebox{1}[0.95]{\includegraphics[width=\textwidth]{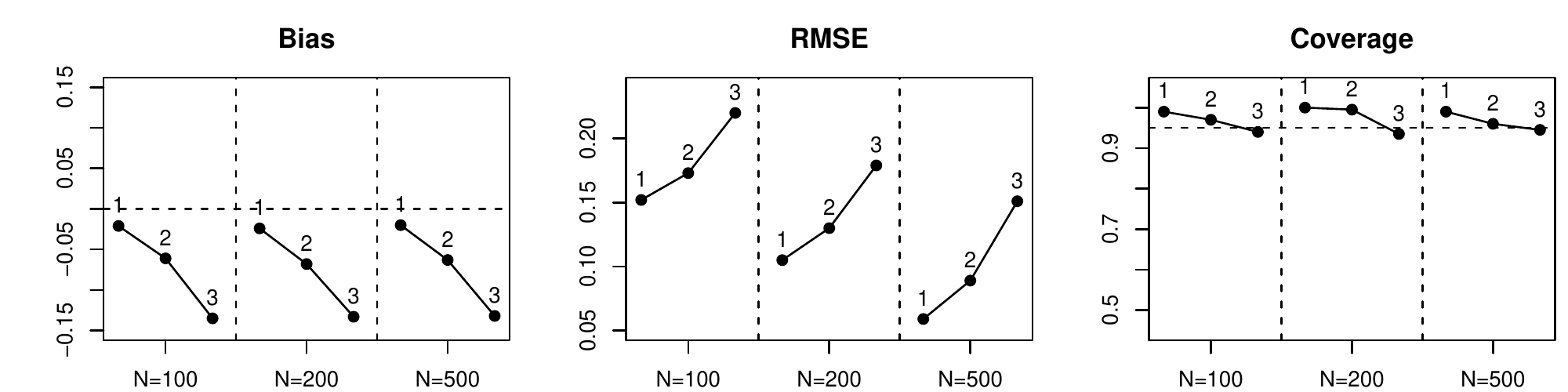}}}
  \caption{Simulations for homogeneity (a) and (b) and heterogeneity (c) and (d). Each subgraph presents the bias, RMSE, coverage proportions of the $95\%$ credible intervals of  $ACE_{\bar{s}\bar{s}}$. Nine combinations of $N_R$ and $N$ are shown in (a) and (b); and nine combinations of $N$ and $d$ (``1'' for $.01$, ``2'' for $.025$, ``3'' for $.05$) are shown in (c) and (d).}
  \end{figure}

\begin{figure}
  \centering
   \subfigure[Correctly specified model without monotonicity]{
   \label{fig:sim2:mon:LH}
    \includegraphics[width=\textwidth]{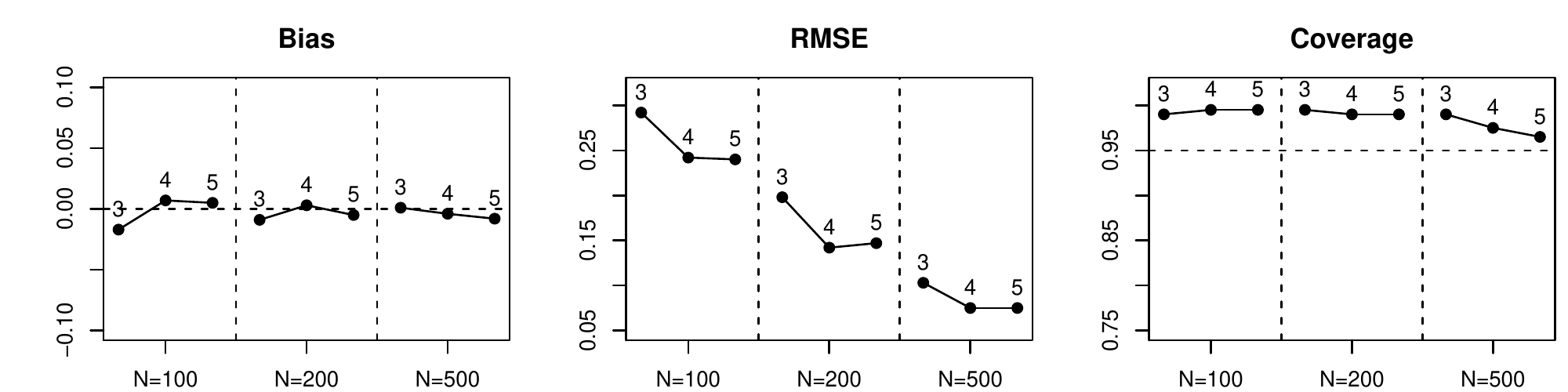}}
    \subfigure[Misspecified model without monotonicity]{
    \label{fig:sim2:non:LH}
  \includegraphics[width=\textwidth]{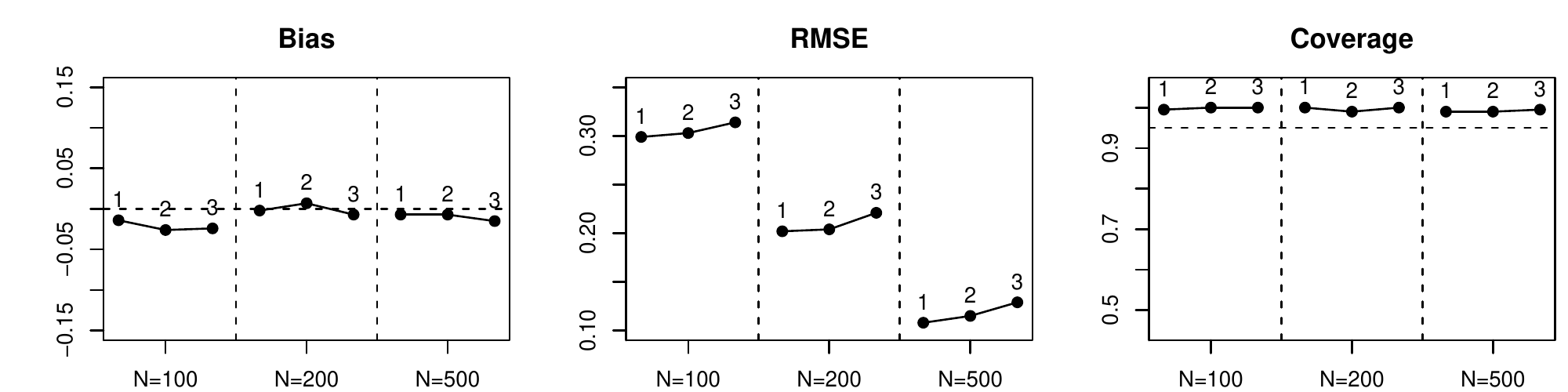}}
  \caption{Simulations for homogeneity (a) and (b) and heterogeneity (c) and (d). Each subgraph presents the bias, RMSE, coverage proportions of the $95\%$ credible intervals of  $ACE_{\bar{s}s}$. Nine combinations of $N_R$ and $N$ are shown in (a) and (b); and nine combinations of $N$ and $d$ (``1'' for $.01$, ``2'' for $.025$, ``3'' for $.05$) are shown in (c) and (d).}
  \end{figure}

 \end{document}